\documentclass[aos,MSNbibl,nameyear,dvips]{arximspdf}

\doi{10.1214/11-AOS883}
\volume{39}
\issue{3}
\pubyear{2011}
\firstpage{1608}
\lastpage{1632}

\makeatletter

\newtheorem{theorem}{Theorem}
\newtheorem{lemma}{Lemma}
\newproclaim{definition}{Definition}
\newproclaim{remarks}{Remarks}

\renewcommand{\kappa}{\varkappa}
\renewcommand{\epsilon}{\varepsilon}
\newcommand{\rd}{\mathrm{d}}
\newcommand{\cA}{\mathcal{A}}
\newcommand{\cB}{\mathcal{B}}
\newcommand{\cF}{\mathcal{F}}
\newcommand{\cH}{\mathcal{H}}
\newcommand{\cL}{\mathcal{L}}
\newcommand{\cR}{\mathcal{R}}
\newcommand{\cW}{\mathcal{W}}
\newcommand{\cX}{\mathcal{X}}

\newcommand{\citeasnoun}[1]{\citet{#1}}

\newcommand{\bB}{\mathbb B}
\newcommand{\bE}{\mathbb E}
\newcommand{\bF}{\mathbb F}
\newcommand{\bH}{\mathbb H}
\newcommand{\bI}{{\mathbb I}}
\newcommand{\bL}{\mathbb{L}}
\newcommand{\bP}{{\mathbb P}}
\newcommand{\bR}{\mathbb{R}}
\newcommand{\mL}{\mathfrak{L}}

\makeatother

\begin{document}
\begin{frontmatter}

\title{Bandwidth selection in kernel density estimation:
Oracle inequalities and adaptive minimax optimality}
\runtitle{Selection of kernel density estimators}

\begin{aug}
\author[A]{\fnms{Alexander} \snm{Goldenshluger}\corref{}\thanksref{t1}\ead[label=e1]{goldensh@stat.haifa.ac.il}}
and
\author[B]{\fnms{Oleg} \snm{Lepski}\ead[label=e2]{lepski@cmi.univ-mrs.fr}}
\thankstext{t1}{Supported by ISF Grant 389/07.}
\runauthor{A. Goldenshluger and O. Lepski}
\affiliation{University of Haifa and Universit\'{e} de Provence}
\address[A]{Department of Statistics
\\
University of Haifa
\\
31905 Haifa\\
Israel
\\
\printead{e1}} 
\address[B]{Laboratoire
d'Analyse, Topologie\\
\quad et Probabilit\'{e}s\\
UMR CNRS 6632
\\
Universit\'{e} de Provence\\
39, rue F. Joliot Curie
\\
13453 Marseille\\
France\\
\printead{e2}}
\end{aug}

\received{\smonth{9} \syear{2010}}
\revised{\smonth{1} \syear{2011}}

%
\begin{abstract}
We address the problem of density estimation with $\mathbb{L}_s$-loss
by selection
of kernel estimators. We develop a selection procedure and derive
corresponding $\mathbb{L}_s$-risk oracle inequalities.
It is shown
that the proposed selection rule leads to the
estimator being minimax adaptive over a scale of the anisotropic Nikol'skii
classes.
The main technical tools used in our
derivations are uniform bounds
on the $\mathbb{L}_s$-norms of empirical processes developed recently by
Goldenshluger and Lepski [\textit{Ann. Probab.} (2011), to appear].
\end{abstract}

%

\end{frontmatter}

\section{Introduction}\label{secintroduction}
Let $X$ be a random variable in $\bR^d$
having density $f$ with respect to the Lebesgue measure.
We want to estimate $f$ on the basis of the i.i.d. sample
$\cX_n=(X_1,\ldots, X_n)$ drawn from $f$.
Any $\cX_n$-measurable map $\hat{f}\dvtx\bR^d\to\bL_s (\bR^d )$ is
understood as an estimator of $f$,
and its accuracy is measured by the
$\bL_s$-risk:
\[
\cR_s[\hat{f}, f] := [\bE_f\|\hat{f} - f\|_s^q ]^{1/q},\qquad s\in[1,
\infty), q\geq1,
\]
where $\bE_f$ is the expectation with respect to the probability measure
$\bP_f$
of the observations $\cX_n$.
The objective is to develop an estimator of $f$ with small $\bL_s$-risk.

Kernel density estimates originate in
\citeasnoun{rosenblatt} and \citeasnoun{parzen}; this
is one of the most popular techniques
for estimating densities [\citeasnoun{silverman},
\citeasnoun{devroye-gyorfi}].
Let $K\dvtx\bR^d\to\bR$ be a fixed function such that $\int
K(x)\,\rd x=1$ (we
call such functions \textit{kernels}).
Given a \textit{bandwidth vector} $h=(h_1,\ldots, h_d)$, $h_i>0$,
the kernel estimator $\hat{f}_h$ of $f$
is defined by
%
\begin{equation}\label{eqkernel-est}
\hat{f}_h (t)= \frac{1}{nV_h} \sum_{i=1}^n K \biggl(\frac{t-X_i}{h} \biggr)=
\frac{1}{n}\sum_{i=1}^n K_h(t-X_i),
\end{equation}
where
$V_h:=\prod_{i=1}^d h_i$,
$u/v$ for $u,v\in\bR^d$ stands for the coordinate-wise division,
and
$K_h(\cdot):=V_h^{-1} K(\cdot/h)$.
It is well known that accuracy properties of $\hat{f}_h$ are
determined by the
choice of the bandwidth $h$,
and \textit{bandwidth selection} is the central
problem in kernel density estimation.
There are different approaches to the problem
of bandwidth selection.

The \textit{minimax approach} is based on
the assumption that
$f$ belongs to a given class
of densities $\bF$, and
accuracy of $\hat{f}_h$ is measured by
its maximal $\bL_s$-risk over the class $\bF$,
\[
\cR_s[\hat{f}_h; \bF] := \sup_{f\in\bF} \cR_s[\hat{f}_h; f].
\]
Typically $\bF$ is a class of smooth functions, for example,
the H\"{o}lder, Nikol'skii or Besov functional class. Then the
bandwidth $h$ is
selected so that the maximal risk $\cR_s[\hat{f}_h; \bF]$
(or a reasonable upper bound on it) is minimized with respect to $h$.
Such a choice leads to a
deterministic bandwidth $h$ depending on the sample size $n$,
and on the underlying functional class $\bF$.
In many cases the resulting kernel estimator constructed in this way
is \textit{rate optimal} (or \textit{optimal in order}) over the class $\bF$.
Minimax kernel density estimation
with $\bL_s$-risks on
$\bR^d$ was considered in
\citeasnoun{bretagnolle}, \citeasnoun{Ibr-Has1}, \citet{Ibr-Has2},
\citeasnoun{devroye-gyorfi}, \citeasnoun{Has-Ibr},
\citeasnoun{Donoho}, \citeasnoun{kerk},
\citeasnoun{Juditsky} and \citeasnoun{mason} where further references
can be found.

The \textit{oracle approach} considers a set of kernel estimators
$\cF(\cH)=\{\hat{f}_h, h\in\cH\}$,
and aims at a measurable data-driven choice $\hat{h}\in\cH$ such
that for every $f$ from a \textit{large} functional class
the following {\em$\bL_s$-risk oracle inequality} holds:
%
\begin{equation}\label{eqoracle}
\cR_s[\hat{f}_{\hat{h}}; f] \leq C \inf_{h\in\cH}\cR_s[\hat
{f}_h;f] +
\delta_n.
\end{equation}
Here $C>0$ is a constant independent of $f$ and $n$, and
the remainder $\delta_n$ does not depend on $f$.
Oracle inequalities with ``small'' remainder term
$\delta_n$ and constant $C$ close to 1 are of prime interest; they are
key tools for
establishing minimax and adaptive minimax results in estimation problems.
To the best of our knowledge, oracle
inequalities of the type (\ref{eqoracle})
were
established only in the cases $s=1$ and \mbox{$s=2$}.
\citeauthor{dev-lug96} (\citeyear{dev-lug96}, \citeyear{dev-lug97}, \citeyear{devroye-lugosi}) established oracle inequalities for $s=1$.
The case \mbox{$s=2$} was studied by
\citeauthor{massart} [\citeyear{massart}, Chapter~7],
\citeasnoun{samarov},
\citeasnoun{rigollet-tsybakov} and \citeasnoun{birge}. The last cited paper
contains a detailed discussion of recent developments in this area.

The contribution of this paper is twofold.
First, we propose a selection procedure for a set of kernel
estimators, and establish for the corresponding $\bL_s$-risk, $s\in
[1,\infty)$,
oracle inequalities
of the type (\ref{eqoracle}).
Second, we demonstrate that
our selection rule
leads to a minimax adaptive estimator over a scale of
the anisotropic Nikol'skii classes (see Section~\ref{secnikolski}
below for the
class definition).

More specifically, let $h^{\min}= (h^{\min}_{1},\ldots,h^{\min}_{d}
)$ and
$h^{\max}= (h^{\max}_{1},\ldots, h^{\max}_{d} )$ be two fixed
vectors satisfying
$0<h^{\min}_{i}\leq h^{\max}_{i}\leq1$, $\forall i$, and let
%
\begin{equation}\label{eqcH}
\cH:=\bigotimes_{i=1}^d [h^{\min}_i, h^{\max}_i ].
\end{equation}
Consider the set of kernel estimators
%
\begin{equation}\label{eqcF}
\cF(\cH)=\{\hat{f}_h, h\in\cH\},
\end{equation}
where $\hat{f}_h$ is given in (\ref{eqkernel-est}).
We propose a measurable choice $\hat{h}\in\cH$
such that the resulting estimator
$\hat{f}=\hat{f}_{\hat{h}}$ satisfies the following oracle inequality:
%
\begin{equation}\label{eqoracle-simple}
\cR_s[\hat{f}_{\hat{h}}; f] \leq
\inf_{h\in\cH} \{ (1+3\|K\|_1 )\cR_s [\hat{f}_h; f] +
C_s(nV_h)^{-\gamma_s} \} + \delta_{n,s}.
\end{equation}
The constants $C_s$, $\gamma_s$, and the remainder term $\delta_{n,s}$
admit different expressions
depending on the value of $s$.
\begin{itemize}
\item If $s\in[1,2)$, then (\ref{eqoracle-simple}) holds for all
densities $f$
with $\gamma_s=1-\frac{1}{s}$, $C_s$
depending on the kernel $K$ only, and with
\[
\delta_{n,s}=c_1 (\ln n)^{c_2}n^{1/s}\exp\{-c_3 n^{2/s-1} \}
\]
for some constants
$c_i$, $i=1,2,3$.
\item If $s\in[2,\infty)$, then (\ref{eqoracle-simple}) holds for
all densities
$f$ uniformly bounded by a constant $\mathrm{f}_\infty$
with $\gamma_s=\frac{1}{2}$, $C_s$ depending on $K$ and $\mathrm
{f}_\infty$
only, and with
\[
\delta_{n,s}=c_1 (\ln n)^{c_2}n^{1/2}\exp\{-c_3 V_{\max}^{-2/s} \},\qquad
V_{\max}:=V_{h^{\max}},
\]
for some constants
$c_i$, $i=1,2,3$.
We emphasize that the proposed
selection rule is fully data-driven and does not use information on the
value of
$\mathrm{f}_\infty$.
\end{itemize}
Thus,
the oracle inequality (\ref{eqoracle-simple}) holds with exponentially small
(in terms of dependence on $n$)
remainder $\delta_{n,s}$
(by choice of $V_{\max}$ in the case $s\in[2,\infty)$).
We stress that explicit nonasymptotic expressions for
$C_s$, $c_1$, $c_2$ and $c_3$ are available.
It is important to realize that the term $C_s(nV_h)^{-\gamma_s}$ is a tight
upper bound on the stochastic error of the kernel estimator $\hat
{f}_h$. This
fact allows to derive rate optimal estimators that
adapt to unknown smoothness of the density $f$.
In particular, in Section~\ref{secnikolski}
we apply our oracle
inequalities in order to develop
a rate optimal adaptive kernel estimator
for the anisotropic Nikol'skii classes.
Minimax estimation of densities from such classes was studied in
\citeasnoun{Ibr-Has2},
while the problem of adaptive estimation was not considered in the
literature.\vadjust{\goodbreak}

The paper is structured as follows.
In Section~\ref{secrule},
we define our selection rule and prove key oracle inequalities.
Section~\ref{secnikolski} discusses adaptive rate optimal
estimation of
densities for a scale of anisotropic Nikol'skii classes.
Proofs of all results are given in Section~\ref{secproofs}.
%
\section{Selection rule and oracle inequalities}\label{secrule}
Let $\cF(\cH)$ be the set of
kernel density estimators defined in (\ref{eqcF}).
We want to select an estimator from the family~$\cF(\cH)$.
For this purpose, we need to impose some assumptions
and establish notation that will be used in the definition of our
selection procedure.
%
%
\subsection{Assumptions}
The following assumptions on the kernel $K$ will be used
throughout the paper.
\begin{longlist}[(K2)]
\item[(K1)] The kernel $K$ satisfies the Lipschitz condition
\[
|K(x)-K(y)|\leq L_K |x-y|\qquad  \forall x,y\in\bR^d,
\]
where $|\cdot|$ denotes the Euclidean distance. Moreover, $K$ is compactly
supported,
and, without loss of generality, $\operatorname{supp}(K)\subseteq
[-1/2, 1/2]^d$.
\item[(K2)]
There exists a real number
$\mathrm{k}_\infty<\infty$
such that
$\|K\|_\infty\leq\mathrm{k}_\infty$.
\end{longlist}

Assumptions (K1) and (K2) are rather standard in kernel density estimation.
We note that
Assumption (K1) can be weakened in several ways. For example,
it suffices to assume that $K$ belongs
to the isotropic H\"{o}lder ball of functions
$\bH_d(\alpha,L_K)$
with any $\alpha>0$ [in Assumption (K1) $\alpha=1$].

Sometimes we will suppose that
$f\in\bF$, where
\[
\bF:= \biggl\{p\dvtx \bR^d\to\bR\dvtx p\geq0, \int p =1, \|p\|_\infty\leq
\mathrm{f}_\infty<\infty\biggr\},
\]
and $\mathrm{f}_\infty$ is a fixed constant.
Without loss of generality we assume that $\mathrm{f}_\infty\geq1$.
%
\subsection{Notation}
For any $U\dvtx\bR^d\to\bR$ and $s\in[1,\infty)$
define
\[
\rho_s(U):= \cases{
4n^{1/s-1} \|U\|_s, &\quad $ s\in[1,2)$,\cr
n^{-1/2} \|U\|_2, &\quad $ s=2$,
}
\]
and if $s\in(2,\infty)$, then we set
\[
\rho_s(U) := D_s \biggl\{n^{-1/2} \biggl(\int\biggl[
\int U^2(t-x) f(x) \,\rd x \biggr]^{s/2} \,\rd t \biggr)^{1/s} + 2
n^{1/s-1}\|U\|_s \biggr\},
\]
where $D_s:=15s/\ln s$ is the best-known constant in the Rosenthal inequality
[\citet{Johnson}].
Observe that $\rho_s(U)$ depends on $f$ when $s\in(2,\infty)$; hence
we will
also consider the empirical counterpart of $\rho_s(U)$:
\[
\hat{\rho}_s(U):= D_s \Biggl\{n^{-1/2} \Biggl(\int\Biggl[\frac{1}{n}\sum_{i=1}^n
U^2(t-X_i) \Biggr]^{s/2} \,\rd t \Biggr)^{1/s} + 2 n^{1/s-1}\|U\|_s \Biggr\}.\vadjust{\goodbreak}
\]
We put also
\[
r_s(U) := \rho_s(U) \vee n^{-1/2}\|U\|_2,\qquad
\hat{r}_s(U) := \hat{\rho}_s(U) \vee n^{-1/2}\|U\|_2
\]
and
\[
g_s(U):= \cases{
32 \rho_s(U), &\quad $ s\in[1,2)$,\vspace*{2pt}\cr
\displaystyle \frac{25}{3}\rho_2(U), &\quad $ s=2$,\vspace*{2pt}\cr
32 \hat{r}_s(U), &\quad $ s>2$.
}
\]

Armed with this notation we are ready to describe our selection rule.

\subsection{Selection rule}
The rule is based on auxiliary
estimators $\{\hat{f}_{h,\eta}, h, \eta\in\cH\}$
that are defined as follows:
for every pair $h, \eta\in\cH$ we let
\[
\hat{f}_{h,\eta}(t):=
\frac{1}{n}\sum_{i=1}^n [K_h*K_\eta](t-X_i),
\]
where ``$*$'' stands for the convolution on $\bR^d$.
Define also
\begin{eqnarray}\label{eqmajor}
 m_s(h,\eta)&:=&g_s(K_\eta) + g_s(K_h*K_\eta)\qquad  \forall h,
\eta\in\cH,
\nonumber\\[-8pt]\\[-8pt]
m_s^*(h)&:=& \sup_{\eta\in\cH} m_s(\eta, h)\qquad  \forall h\in\cH.\nonumber
\end{eqnarray}
For
every $h \in\cH$ let
%
\begin{eqnarray}\label{eqrule-0}
\hat{R}_{h} := \sup_{\eta\in\cH}
[ \|\hat{f}_{h, \eta}- \hat{f}_{\eta}\|_s -
m_s(h, \eta) ]_+ + m^*_s(h).
\end{eqnarray}
The selected bandwidth $\hat{h}$ and the corresponding kernel density
estimator
are defined by
%
\begin{equation}\label{eqrule}
\hat{h}:=\operatorname{arg}\inf_{h\in\cH} \hat{R}_{h},\qquad
\hat{f}=\hat{f}_{\hat{h}}.
\end{equation}

The selection rule (\ref{eqmajor})--(\ref{eqrule}) is a
refinement of
the one introduced recently in
\citeauthor{GL-1} (\citeyear{GL-1}, \citeyear{GL-2}) for the Gaussian white noise model.

\begin{remarks*}
1. It is easy to check that Assumption (K1) implies that $\hat{R}_{h}$ and
$m_s^*(h)$ are continuous random functions on the
compact subset $\cH\subset\bR^d$. Thus, $\hat{h}$ exists and is
measurable
[\citet{jennrich}].

2.
We call function $m_s(\cdot, \cdot)$ the \textit{majorant}.
In fact, if $\xi_h$ and $\xi_{h,\eta}$ denote the stochastic errors of
estimators $\hat{f}_h$ and $\hat{f}_{h,\eta}$, respectively, that is, if
\begin{eqnarray*}
\xi_{h} (t) &:=& \frac{1}{n}\sum_{i=1}^n [K_h(t-X_i) - \bE_f
K_h(t-X) ],
\\
\xi_{h,\eta}(t) &:=& \frac{1}{n}
\sum_{i=1}^n \{ [K_h*K_\eta](t-X_i) - \bE_f [K_h*K_\eta](t-X) \},
\nonumber
\end{eqnarray*}
then it is seen from the proofs of Theorems~\ref{ths<=2} and \ref
{ths>2} below
that
$m_s(h,\eta)$ uniformly ``majorates'' $\|\xi_{h,\eta}-\xi_\eta\|
_s$ in the
sense that the expectation
\[
\bE_f \sup_{(h, \eta)\in\cH\times\cH}
[\|\xi_{h, \eta}-\xi_{\eta}\|_s - m_{s}(h, \eta) ]_+^q
\]
is ``small.''

3.
It is important to realize that the majorant $m_s(h,\eta)$ is
explicitly given
and does not depend
on the density $f$ to be estimated. The majorant is completely
determined by
kernel $K$ and observations, and thus it is
available to the statistician.
\end{remarks*}

\subsection{Oracle inequalities}
Now we are in a position to establish oracle inequalities
on the risk of the estimator $\hat{f}=\hat{f}_{\hat{h}}$
given by (\ref{eqrule}). Put
\[
A_\cH:=
\prod_{i=1}^d [1\vee\ln (h_i^{\max}/h_i^{\min} ) ],\qquad
B_\cH:= [1\vee\log_2 (V_{\max}/V_{\min} ) ],
\]
where from now on
\[
V_{\min}:=\prod_{i=1}^d h_i^{\min},\qquad
V_{\max}:=\prod_{i=1}^d h_i^{\max}.
\]

The next two statements, Theorems~\ref{ths<=2} and~\ref{ths>2},
provide oracle inequalities on the $\bL_s$-risk of
$\hat{f}$ in the cases $s\in[1,2]$ and $s\in(2,\infty)$, respectively.
\begin{theorem}\label{ths<=2}
Let Assumptions~\textup{(K1)} and \textup{(K2)} hold.
\begin{itemize}[(ii)]
\item[(i)] If $s\in[1,2)$, then for all $f$ and $n\geq4^{2s/(2-s)}$
\begin{eqnarray}\label{eqor-1}
\cR_s[\hat{f}; f] &\leq& \inf_{h\in\cH} \bigl[(1+3\|K\|_1)\cR_s[\hat{f}_h,
f]+
C_1(nV_h)^{1/s -1} \bigr]
\nonumber\\[-8pt]\\[-8pt]
&&{}
+
C_2 A_\cH^{4/q} n^{1/s} \exp\biggl\{-\frac{2n^{2/s-1}}{37q} \biggr\}.
\nonumber
\end{eqnarray}
\item[(ii)] If $s=2$ and $\mathrm{f}_\infty^2V_{\max}
+4n^{-1/2}\leq1/8$,
then
for all $f\in\bF$
\begin{eqnarray}\label{eqor-2}
\cR_s[\hat{f}; f] &\leq& \inf_{h\in\cH} \bigl[(1+3\|K\|_1)\cR_s[\hat{f}_h,
f]+
C_3(nV_h)^{-1/2} \bigr]
\nonumber\\[-8pt]\\[-8pt]
&&{} +
C_4 A_\cH^{4/q} n^{1/2} \exp\biggl\{-\frac{1}{16q[\mathrm{f}^2_\infty
V_{\max} +4n^{-1/2}]} \biggr\}.\nonumber
\end{eqnarray}
\end{itemize}
Here $C_1$ and $C_3$ are absolute constants, while $C_2$ and $C_4$ depend
on $L_K$, $\mathrm{k}_\infty$, $d$ and $q$ only.
\end{theorem}
%
%
\begin{theorem}\label{ths>2}
Let Assumptions \textup{(K1)} and \textup{(K2)} hold,
$s\in(2,\infty)$, and assume
that for
some
$C_1=C_1(K,s, d)>1$
\[
nV_{\min}> C_1,\qquad  V_{\max}\geq1/\sqrt{n}.\vadjust{\goodbreak}
\]
If $n\geq C_2$ for some constant $C_2$ depending on $L_K$, $\mathrm
{k}_\infty$,
$\mathrm{f}_\infty$, $d$
and $s$ only, then $\forall f\in\bF$,
\begin{eqnarray}\label{eqor-3}
\cR_s[\hat{f};f] &\leq&
\inf_{h\in\cH} \bigl[(1+3\|K\|_1)\cR_s[\hat{f}_h,
f]+
C_3\mathrm{f}_\infty^{1/2}
(nV_h)^{-1/2} \bigr]
\nonumber\\[-8pt]\\[-8pt]
&&{}
+
C_{4} A_\cH^{4/q}B_\cH^{1/q} n^{1/2} [
\exp\{-C_{5}b_{n,s}\} + \exp\{-C_6\mathrm{f}_\infty^{-1} V_{
\max}^{-2/s} \} ],
\nonumber
\end{eqnarray}
where
$b_{n,s}:=n^{4/s-1}$ if $s\in(2,4)$, and
$b_{n,s}:=[\mathrm{f}_\infty V_{\max}^{4/s}]^{-1}$ if $s\geq4$.
The constants $C_i$, $i=3,\ldots,6$, depend on
$L_K$, $\mathrm{k}_\infty$, $d$, $q$ and $s$
only.
\end{theorem}

\begin{remarks*}
1.
All constants appearing in Theorems~\ref{ths<=2} and~\ref{ths>2} can be
expressed explicitly [see 
Lemmas~\ref{props<=2} and~\ref{props>2} below and corresponding results
in Goldenshluger and Lepski (\citeyear{GL}) for details].

2. We will show that for given $h$ the expected value of the
stochastic error
of the estimator $\hat{f}_h$, that is, $ (\bE\|\xi_h\|^{q}_s )^{1/q}$,
admits the upper bound
of the order $O((nV_h)^{1/s-1})$ when $s\in[1,2)$ and $O((nV_h)^{-1/2})$
when $s\in(2,\infty)$. It is also obvious that
\begin{eqnarray*}
\cR_s[\hat{f}_h;f] &\leq& \|B_h\|_s+ (\bE_f\|\xi_h\|^{q}_s )^{1/q},
\end{eqnarray*}
where $B_{h}(f,t) :=\int K_h(t-x) f(x) \,\rd x - f(t), t\in\bR^d$.
Thus, our estimator attains,
up to a constant and remainder term, the minimum of the sum of
the bias and the upper bound on the stochastic error.
This form of the oracle inequality is 
convenient for deriving
minimax and minimax adaptive results (see Section~\ref{secnikolski}).
Indeed, bounds on the bias and
the stochastic error are usually developed separately and
require completely different techniques.

3.
We note that $A_\cH\leq O([\ln n]^d)$
and $B_\cH\leq O(\ln n)$ for any set $\cH\subset[0,1]^d$ such that
$h^{\min}_i\geq
O(n^{-c})$, $c>0$, $\forall i=1,\ldots, d$. If $s\in(2,\infty)$, and
if the set
of considered
bandwidths $\cH$
is such that $V_{\max}=[\kappa\ln n]^{-s/2}$ for some $\kappa>0$,
then the
second term on the right-hand side of (\ref{eqor-2}) and (\ref
{eqor-3}) can
be made negligibly small by carefully choosing the constant $\kappa$.
Observe that conditions ensuring consistency
of $\hat{f}_h$ are $nV_h \to\infty$ and $V_h\to0$ as $n\to\infty$;
thus the requirement $V_{\max}=[\kappa\ln n]^{-s/2}$ is not restrictive.
Note also that
in the case $s\in[1,2)$ the second term on the right-hand side of
(\ref{eqor-1}) is exponentially small in $n$ for any~$\cH$.

4.
The condition $V_{\max}\geq1/\sqrt{n}$ is imposed only
for the sake of convenience in the presentation of our results.
Clearly, we would like to have the set $\cH$ as large as possible; hence
consideration of vectors $h^{\max}$ such that $V_{\max}=V_{h^{\max
}}\leq
1/\sqrt{n}$
does not make much sense.
%

5.
Note that the oracle inequalities (\ref{eqor-1}), (\ref{eqor-2}) and
(\ref{eqor-3}) of
Theorems~\ref{ths<=2} and~\ref{ths>2} hold under very mild conditions
on the
density $f$. In particular, in the case $s\in[1,2)$ the inequality
(\ref{eqor-1}) holds for all densities, and only
boundedness of $f$ is required for (\ref{eqor-2}) and (\ref{eqor-3}).

6. It should be also mentioned that if for $s\in[1,2)$ we impose
additional conditions on $f$ [e.g., such as the domination condition in
\citeasnoun{Donoho}, page~514], then the order of the stochastic
error of $\hat{f}_h$ can be improved to $O((nV_h)^{-1/2})$.
This will lead to the oracle inequality (\ref{eqor-1}) with the term
$C_1(nV_h)^{1/s-1})$ replaced by $C_1(nV_h)^{-1/2}$.
However, $O((nV_h)^{1/s-1})$ is a tight upper bound on the stochastic error
of $\hat{f}_h$ when no conditions on $f$ are assumed.
In particular, it is well known
that smoothness condition alone is not sufficient for consistent density
estimation
on $\bR^d$ with $\bL_1$-losses [\citet{Ibr-Has2}].
\end{remarks*}
%
\subsection{$\bL_s$-risk oracle inequalities}

As it was mentioned above, the oracle inequalities of Theorems~\ref{ths<=2}
and~\ref{ths>2} are useful for derivation of adaptive rate
optimal estimators. 
They are established under very mild assumptions on the density $f$.
However, it is not clear
how the second term under the infimum sign
on the right-hand side of the
developed
oracle inequalities is compared to $\cR_s[\hat{f}_h;f]$.
Traditionally oracle inequalities compare
the risk of a proposed estimator
to the risk of the best estimator in the given family;
cf.~(\ref{eqoracle}).
Therefore the natural question is whether an {\em$\bL_s$-risk oracle
inequality}
of the type (\ref{eqoracle})
can be derived from
the results of Theorems~\ref{ths<=2} and~\ref{ths>2}.

In this section we provide an answer to this question.
We will be mostly interested in finding minimal assumptions 
on the underlying density $f$
that are sufficient for establishing the $\bL_s$-risk
oracle inequality.
It will be shown
that this problem is directly related to establishing a lower bound on
the term
$ (\bE_f\|\xi_h\|^{q}_s )^{1/q}$.

Let
$\mu\in(0,1)$ and $\nu>0$ be fixed real numbers.
Denote by $\bF_{\mu,\nu}$
the set of all probability densities $f$ satisfying the
following condition:
\[
\exists B\in\cB(\bR^d)\dvtx\qquad  \operatorname{mes}(B)\leq\nu,\qquad
\int_{B}f\geq\mu.
\]
Here $\cB(\bR^d)$ is the Borel $\sigma$-algebra on $\bR^d$ and
$\operatorname{mes}(\cdot)$ is the Lebesgue measure on~$\bR^d$.

Below we will assume that $f\in\bF_{\mu,\nu}$ for some $\mu$ and
$\nu$.
This condition is very weak.
For example, if $\cF$ is a set of densities
such that either
(i) $\cF$ is a totally bounded subset of
$\bL_1(\bR^{d})$,
or
(ii) the family of probability measures
$\{\bP_f, f\in\cF\}$ is tight,
then
for any $\mu\in(0,1)$ there exists $0<\nu<\infty$
such that $\cF\subseteq\bF_{\mu,\nu}$.
The statement (i) is a consequence of the Kolmogorov--Riesz compactness theorem.
\begin{theorem}\label{thoracle}
Let $s\in[2,\infty)$ and suppose that assumptions of
Theorems~\ref{ths<=2}\textup{(ii)} and~\ref{ths>2} are fulfilled. If
$s>2$, then
assume additionally
that $f\in\bF_{\mu,\nu}$ for some $\mu$ and $\nu$, and
\[
V_{\max}\leq
2^{-1}\mu\biggl[\frac{\|K\|_2}{\|K\|_1} \biggr]^{2}.
\]
If $n\geq C_1=C_1(L_K, \mathrm{k}_\infty,\mathrm{f}_\infty, d, s)$,
then there exists a constant $C_0>0$ [$C_0=C_0(K)$ if $s=2$ and
$C_0=C_0(K,\mu,\nu,s)$ if $s>2$] such that
\begin{eqnarray*}
\cR_s[\hat{f};f] &\leq& C_0
\inf_{h\in\cH} \cR_s[\hat{f}_{\hat{h}}; f]
\nonumber
\\
&&{}
+
C_{2} A_\cH^{4/q}B_\cH^{1/q} n^{1/2} [
\exp\{-C_{3}b_{n,s}\} + \exp\{-C_4\mathrm{f}_\infty^{-1} V_{
\max}^{-2/s} \} ],
\end{eqnarray*}
where
$b_{n,s}:=n^{4/s-1}$ if $s\in(2,4)$ and
$b_{n,s}:=[\mathrm{f}_\infty V_{\max}^{4/s}]^{-1}$ if $s\geq4$.
The constants $C_i$ depend on
$L_K$, $\mathrm{k}_\infty$, $d$, $q$ and $s$
only.
\end{theorem}

The proof indicates that Theorem~\ref{thoracle}
follows from the fact that for any $s\in[2,\infty)$
one has
%
\begin{equation}\label{eqlow-stoch}
[\bE_f\|\xi_h\|^q_s]^{1/q}\geq c(nV_h)^{-1/2}\qquad  \forall h,
\end{equation}
where $c>0$ is a constant. This lower bound
holds under very weak conditions on the density $f$ (for arbitrary $f$
if $s=2$ and $f\in\bF_{\mu,\nu}$ if $s>2$).
In order to prove the similar $\bL_s$-risk oracle inequality in the case
$s\in[1,2)$ it would be sufficient to show that
$[\bE_f\|\xi_h\|^q_s]^{1/q}\geq c(nV_h)^{-1+1/s}$ for any $h$.
However, the last lower bound cannot hold in such generality as
(\ref{eqlow-stoch}). In particular, according to Remark 5 after
Theorem~\ref{ths>2}, $[\bE_f\|\xi_h\|^q_s]^{1/q}\leq c(nV_h)^{-1/2}$
for all $h$ under a tail
domination condition (e.g., for compactly supported densities).
Under such a domination condition
the corresponding $\bL_s$-risk oracle inequality can be easily established
using the same arguments as in the proof of Theorem~\ref{thoracle}.
%
\subsection{Generalization}
\label{secdiscussion}

Although in the present paper we focus on the bandwidth selection,
the proposed selection rule can be easily extended
to very general families of linear estimators.

Let $\mL$ be the collection of functions $\cL\dvtx\bR^d\times\bR^d\to
\bR$
such that
\[
\int_{\bR^d}\cL(t,x)\,\rd t=1\qquad  \forall x\in\bR^d.
\]
Consider the following family of estimators generated by 
$\mL$:
\[
\cF(\mL)=\Biggl \{\hat{f}_{\cL} (\cdot)=\frac{1}{n}\sum_{i=1}^n
\cL(\cdot,X_i), \cL\in\mL\Biggr\}.
\]
The objective is to propose the selection rule from the family
$\cF(\mL)$ and to establish for the obtained estimator $\bL_s$-oracle
inequality.
A close
inspection of the proofs of Theorems~\ref{ths<=2} and~\ref{ths>2}
leads to the following generalization
of the selection rule
(\ref{eqrule}).

For any couple $\cL,\cL^\prime\in\mL$ let 
\[
[\cL\otimes\cL^\prime](t,x):=\int_{\bR^d}\cL(t,y)\cL^\prime
(y,x)\,\rd y
\]
and define the estimator
\[
\hat{f}_{\cL\otimes\cL^\prime} (\cdot)=\frac{1}{n}\sum_{i=1}^n
[\cL\otimes\cL^\prime](\cdot,X_i).
\]
Let
\begin{eqnarray*}
\xi_{\cL} (t) &:=& \frac{1}{n}\sum_{i=1}^n [\cL(t,X_i) - \bE_f
\cL(t,X) ],
\\
\xi_{\cL\otimes\cL^\prime}(t) &:=& \frac{1}{n}
\sum_{i=1}^n \{ [\cL\otimes\cL^\prime](t,X_i) - \bE_f
[\cL\otimes\cL^\prime](t,X) \}.
\nonumber
\end{eqnarray*}
Suppose that for any $\cL,\cL^\prime\in\mL$ one can find
a majorant
$m_s (\cL,\cL^\prime)$
for $\|\xi_{\cL\otimes\cL^\prime}-\xi_{\cL^\prime}\|_s$. In
other words,
suppose that
the expectation
\[
\bE_f \sup_{(\cL, \cL^\prime)\in\mL\times\mL}
[\|\xi_{\cL\otimes\cL^\prime}-\xi_{\cL^\prime}\|_s -m_s (\cL
,\cL^\prime)
]_+^q
\]
is ``small,'' and
analogues of
Lemmas~\ref{props<=2} and~\ref{props>2} given
below are proved.
We refer
to Goldenshluger and Lepski (\citeyear{GL}), where
results of this type for various collections $\mL$
can be found.
%

For
every $\cL\in\mL$ let
%
\begin{eqnarray}\label{eqrule-0-general}
&& \hat{R}_{\cL} := \sup_{\cL^{\prime}\in\mL}
[ \|\hat{f}_{\cL\otimes\cL^\prime}- \hat{f}_{\cL^\prime}\|_s -
m_s(\cL,\cL^\prime) ]_+ + \sup_{\cL^\prime\in\mL} m_s(\cL
^\prime, \cL) ,
\end{eqnarray}
and define
%
\begin{equation}\label{eqrule-general}
\widehat{\cL}:=\operatorname{arg}\inf_{\cL\in\mL} \hat{R}_{\cL}.
\end{equation}
The selected estimator is $\hat{f}=\hat{f}_{\widehat{\cL}}$.

In order to prove
analogues of Theorems~\ref{ths<=2}
and~\ref{ths>2} the
following assumption ({\it commutativity property}) on the
collection $\mL$ has to be imposed:
%
\begin{equation}
\label{eqcommutativity-property}
\int_{\bR^d}\cL(\cdot,y)\cL^\prime(y,\cdot)\,\rd y =
\int_{\bR^d}\cL^\prime(\cdot,y)\cL(y,\cdot)\,\rd y\qquad  \forall
\cL,\cL^\prime\in\mL.
\end{equation}
Thus, using
the commutativity property (\ref{eqcommutativity-property})
and majorants for the $\bL_s$-norms of empirical processes derived
in Goldenshluger and Lepski (\citeyear{GL}), one can
establish $\bL_s$-oracle inequalities for the selection rule
(\ref{eqrule-0-general})--(\ref{eqrule-general}).
%
\section{Adaptive estimation of densities with anisotropic smoothness}
\label{secnikolski}
In this section we illustrate the use of oracle inequalities
of Theorems~\ref{ths<=2} and~\ref{ths>2} for derivation
of adaptive rate optimal density estimators.

We start with the definition of the \textit{anisotropic
Nikol'skii class of functions}.
\begin{definition}
Let $p\in[1,\infty]$, $\alpha=(\alpha_1,\ldots, \alpha_d)$,
$\alpha_i>0$, and
$L>0$.
We say that a density $f\dvtx\bR^d\to\bR$ belongs to the
anisotropic Nikol'skii class $N_{p, d}(\alpha, L)$ of functions
if:\vspace*{-2pt}
\begin{longlist}[(ii)]
\item[(i)]
$\|D_i^{\lfloor\alpha_i\rfloor} f\|_p \leq L$,
for all $i=1,\ldots, d$;

\item[(ii)] for all $i=1,\ldots, d$, and all $z\in\bR^1$
\begin{eqnarray*}
&&\biggl\{
\int \bigl|D_i^{\lfloor\alpha_i\rfloor}f(t_1,\ldots,t_i+z,\ldots,t_d)
- D_i^{\lfloor\alpha_i\rfloor}f(t_1,\ldots,t_i,\ldots,t_d) \bigr|^p
\,\rd t
\biggr\}^{1/p}
\\[-2pt]
&&\qquad \leq L|z|^{\alpha_i-\lfloor\alpha_i\rfloor}.
\end{eqnarray*}
\end{longlist}
Here $D_i^kf$ denotes the $k$th-order partial derivative of $f$ with
respect to
the variable $t_i$ and $\lfloor\alpha_i\rfloor$ is the largest
integer strictly
less than $\alpha_i$.\vspace*{-2pt}
\end{definition}

The functional classes $N_{p,d}(\alpha, L)$ were considered in approximation
theory by
Nikol'skii; see, for example, \citeasnoun{nikol}.
Minimax
estimation of densities from the class $N_{p,d}(\alpha, L)$
was considered in
\citeasnoun{Ibr-Has2}. We refer also to
\citeasnoun{lepski-kerk} where the problem of adaptive estimation over a
scale of classes
$N_{p,d}(\alpha,L)$ was treated for the Gaussian white noise model.

Consider the following family of kernel estimators.
Let $u$ be
an integrable, compactly supported
function on $\bR$ such that $\int u(y)\, \rd y=1$. As in \citeasnoun
{lepski-kerk},
for some integer number $l$ we put
\[
u_l(y) := \sum_{k=1}^l \pmatrix{l\cr k}
(-1)^{k+1}
\frac{1}{k}u \biggl(\frac{y}{k} \biggr),\vspace*{-2pt}
\]
and define
%
\begin{equation}\label{eqkernel-K}
K(t):=\prod_{i=1}^d u_l(t_i),\qquad  t=(t_1,\ldots, t_d).
\end{equation}
The kernel $K$ constructed in this way is bounded and compactly
supported, and
it is easily verified that
\[
\int K(t) \,\rd t=1,\qquad  \int K(t) t^k \,\rd t=0\qquad  \forall|k|=1,\ldots, l-1,
\]
where $k=(k_1,\ldots, k_d)$ is the multi-index, $k_i\geq0$,
$|k|=k_1+\cdots+k_d$ and $t^k= t_1^{k_1}\cdots t_d^{k_d}$ for
$t=(t_1,\ldots, t_d)$.

For fixed $\alpha=(\alpha_1,\ldots, \alpha_d)$ set
$1/\bar{\alpha}=\sum_{i=1}^d (1/\alpha_i)$ and define
\[
\varphi_{n,s}(\bar{\alpha}) :=L^{-\gamma_s/(\bar{\alpha}+\gamma_s)}
n^{-\gamma_s
\bar{\alpha}/(\bar{\alpha}+\gamma_s)},\qquad
\gamma_s:= \cases{
1-1/s, &\quad $ s\in(1,2]$,\cr
1/2, & \quad $s\in(2,\infty)$.
}\vspace*{-2pt}
\]
\begin{theorem}\label{thnikolski}
Let $\cF(\cH)$ be the family of kernel estimators defined
in (\ref{eqkernel-est}), (\ref{eqcH}) and (\ref{eqcF}) that is associated
with the kernel
(\ref{eqkernel-K}).
Let $\hat{f}$ denote the estimator
given by selection according to our rule (\ref{eqmajor})--(\ref{eqrule})
from the family~$\cF(\cH)$.\vadjust{\goodbreak}
\begin{longlist}[(ii)]
\item[(i)] Let $s\in(1,2)$, and assume that $h_i^{\min}=1/n$ and
$h^{\max}_i=1$, $\forall i=1,\ldots, d$.
Then for any class $N_{s,d}(\alpha, L)$ such that $\max_{i=1,\ldots,d}
\lfloor\alpha_i\rfloor\leq l-1$, $L>0$ one has
\[
\limsup_{n\to\infty} \{[\varphi_{n,s}(\bar{\alpha})]^{-1} \cR
_s[\hat{f};
N_{s,d}(\alpha, L)] \} <\infty.
\]
\item[(ii)] Let $s\in[2,\infty)$, and assume that
$h_i^{\min}=\kappa_1/n$ and $h^{\max}_i= [\kappa_2 \ln
n]^{-s/(2d)}$, $\forall
i=1,\ldots, d$ for some constants $\kappa_1$ and $\kappa_2$.
Then for any class $N_{s,d}(\alpha, L)$ such that $\max_{i=1,\ldots,d}
\lfloor\alpha_i\rfloor\leq l-1$, $L>0$ one has
\[
\limsup_{n\to\infty} \{[\varphi_{n,s}(\bar{\alpha})]^{-1} \cR
_s[\hat{f};
N_{s,d}(\alpha, L)] \} <\infty.
\]
\end{longlist}
\end{theorem}

It is well known that $\varphi_{n,s}(\bar{\alpha})$
is the minimax rate of convergence in estimation of densities from
the class $N_{s,d}(\alpha, L)$ [see \citeasnoun{Ibr-Has2} and
\citeasnoun{Has-Ibr}].
Therefore Theorem~\ref{thnikolski} shows that our estimator
$\hat{f}$ is adaptive minimax over a scale of the classes
$N_{s,d}(\alpha, L)$ indexed by $\alpha$ and $L$.

The above result
holds when both the smoothness
and the accuracy are measured in the same
$\bL_s$-norm.
We demonstrate below that if
the additional condition of compact support is imposed,
then the resulting estimator
is adaptive minimax over a much larger scale of
functional classes.
\begin{definition}
Let $p\in[1,\infty]$, $\alpha=(\alpha_1,\ldots, \alpha_d)$,
$\alpha_i>0$,
$L>0$, and let $Q$ be a fixed cube in $\bR^d$.
We say that a density $f\dvtx\bR^d\to\bR$ belongs to the
functional class $W_{p, d}(\alpha, L, Q)$ if
$f\in N_{p,d}(\alpha, L)$, and $\operatorname{supp}(f)\subseteq Q$.
\end{definition}
\begin{theorem}\label{thcompact}
Let $s\in[1,\infty)$, and assume that
$h_i^{\min}=\kappa_1/n$ and $h^{\max}_i= [\kappa_2 \ln n]^{-[s\vee
2]/(2d)}$, $\forall
i=1,\ldots, d$ for some constants $\kappa_1$ and $\kappa_2$.
Let $\cF(\cH)$ be the corresponding family of kernel estimators
that is associated
with the kernel
(\ref{eqkernel-K}).
Let $\hat{f}$ denote the estimator
given by the selection procedure
(\ref{eqmajor})--(\ref{eqrule})
with $s$ substituted by $s\vee2$.
%
%
%
Then for any class $W_{p,d}(\alpha, L, Q)$ such that $p\geq [s\vee2]$,
$\max_{i=1,\ldots,d}
\lfloor\alpha_i\rfloor\leq l-1$, $L>0$
%
\[
\limsup_{n\to\infty} \{[\psi_{n,s}(\bar{\alpha})]^{-1} \cR_s
[\hat{f};
W_{p,d}(\alpha, L, Q) ] \} <\infty,
\]
where
\[
\psi_{n,s}(\bar{\alpha}):= \bigl(L [\operatorname{mes}\{Q\} ]^{(p-[s\vee
2])/p[s\vee2]} \bigr)^{1/(2\bar{\alpha}+1)}
n^{-\bar{\alpha}/(2\bar{\alpha}+1)} .
\]
\end{theorem}

Theorem~\ref{thcompact} shows that if $s\in[1,\infty)$, then the estimator
$\hat{f}$ given by our selection procedure achieves the minimax rate of
convergence
simultaneously on
every class $W_{p,d}(\alpha, L, Q)$ with any $p\geq[s\vee2]$, $\max
_{i=1,\ldots,d}
\lfloor\alpha_i\rfloor\leq l-1$, $L>0$ and any fixed support $Q$. It
should be especially
stressed that no information about the support set $Q$ and the index
$p$ are used in
construction of $\hat{f}$.
%
\section{Proofs}
\label{secproofs}

First we recall that the
accuracy of estimators $\hat{f}_{h}$ and $\hat{f}_{h,\eta}$,
$h, \eta\in\cH$, is characterized by
the bias and stochastic error given by
\begin{eqnarray*}
B_{h}(f, t) &:=& \int K_h(t-x) f(x) \,\rd x - f(t),
\\[-2pt]
\xi_{h} (t) &:=& \frac{1}{n}\sum_{i=1}^n [K_h(t-X_i) - \bE_f
K_h(t-X) ]
\end{eqnarray*}
and
\begin{eqnarray*}
B_{h, \eta}(f, t)&:=& \int[K_h*K_\eta](t-x) f(x) \,\rd x - f(t),
\\[-2pt]
\xi_{h,\eta}(t) &:=& \frac{1}{n}
\sum_{i=1}^n \{ [K_h*K_\eta](t-X_i) - \bE_f [K_h*K_\eta](t-X) \},
\end{eqnarray*}
respectively.

The proofs extensively use results from Goldenshluger and Lepski
(\citeyear{GL}); in what
follows
for the sake of brevity
we refer to this paper as GL (\citeyear{GL}).

\subsection{Auxiliary results}
We start with two auxiliary lemmas that establish probability and
moment bounds
on
$\bL_s$-norms of the processes $\xi_h$ and $\xi_{h,\eta}$.
Proofs of these results are given in the \hyperref[app]{Appendix}.
\begin{lemma}\label{props<=2}
Let Assumptions \textup{(K1)} and \textup{(K2)} hold.
\begin{longlist}[(ii)]
\item[(i)] If $s\in[1,2)$,
then for all $n\geq4^{2s/(2-s)}$ one has
\begin{eqnarray}\label{eqdelta-1-s<2}
&&\Bigl\{\bE_f \sup_{h\in
\cH} [\|\xi_{h}\|_s-32\rho_s(K_h) ]_+^q \Bigr\}^{1/q}\nonumber\\[-9pt]\\[-9pt]
&&\qquad \leq \delta_{n,s}^{(1)}:=
C_1 A^{2/q}_\cH n^{1/s} \exp\biggl\{-\frac{2n^{2/s-1}}{37q} \biggr\},
\nonumber
\\[-2pt]
\label{eqdelta-2-s<2}
&&\Bigl\{\bE_f \sup_{(h,\eta)\in\cH\times
\cH} [\|\xi_{h,\eta}\|_s-32\rho_s(K_h*K_\eta) ]_+^q \Bigr\}^{1/q}\nonumber\\[-9pt]\\[-9pt]
&&\qquad \leq \delta_{n,s}^{(2)}:=
C_2 A_\cH^{4/q} n^{1/s} \exp\biggl\{-\frac{2n^{2/s-1}}{37q} \biggr\}.
\nonumber
\end{eqnarray}
\item[(ii)] Let $f\in\bF$, and assume that
$8[\mathrm{f}_\infty^2V_{\max} +4n^{-1/2}]\leq1$;
then for all $f\in\bF$ one has
\begin{eqnarray}\label{eqdelta-2-1}
&& \biggl\{\bE_f \sup_{h\in\cH} \biggl[\|\xi_{h}\|_2- \frac{25}{3}
\rho_2(K_h) \biggr]_+^q \biggr\}^{1/q}\nonumber\\[-8pt]\\[-8pt]
&&\qquad \leq \delta_{n,2}^{(1)}
:=
C_3 A^{2/q}_\cH n^{1/2} \exp\biggl\{-\frac{1}{16q[V_{\max}\mathrm
{f}^2_\infty
+4n^{-1/2}]} \biggr\},\nonumber
\\\label{eqdelta-2-2}
&&
\biggl\{\bE_f \sup_{(h,\eta)\in\cH\times\cH} \biggl[\|\xi_{h,\eta}\|_2 -
\frac{25}{3} \rho_2(K_h*K_\eta) \biggr]_+^q \biggr\}^{1/q}\nonumber\\[-8pt]\\[-8pt]
&&\qquad \leq
\delta_{n,2}^{(2)}
:=
C_4 A_\cH^{4/q} n^{1/2} \exp\biggl\{-\frac{1}{16q[\mathrm{f}^2_\infty
V_{\max} +4n^{-1/2}]} \biggr\}.\nonumber
\end{eqnarray}
\end{longlist}
The constants $C_i$, $i=1,\ldots,4$,
depend on $L_K$, $\mathrm{k}_\infty$, $d$ and $q$ only.
\end{lemma}

\begin{lemma}\label{props>2}
Let Assumptions \textup{(K1)} and \textup{(K2)} hold, $f\in\bF$, $s>2$, and assume that
\[
n\geq C_1,\qquad  nV_{\min}> C_2
,\qquad  V_{\max}\geq1/\sqrt{n}.
\]
Then the following
statements hold:
\begin{eqnarray}\label{eqlevel-1}
&& \Bigl\{
\bE_f
\sup_{h\in\cH} [\|\xi_{h}\|_s - 32
\hat{r}_s(K_h) ]_+^q \Bigr\}^{1/q}\nonumber\\[-8pt]\\[-8pt]
&&\qquad \leq \delta_{n,s}^{(1)}
:=
C_3 A_\cH^{2/q} B^{1/q}_\cH n^{1/2} \exp\biggl\{-
\frac{C_4}{\mathrm{f}_\infty V_{\max}^{2/s}} \biggr\},
\nonumber
\\\label{eqlevel-2}
&& \Bigl\{\bE_f
\sup_{(h,\eta)\in\cH\times\cH} [\|\xi_{h,\eta}\|_s - 32
\hat{r}_s(K_h*K_\eta) ]_+^q \Bigr\}^{1/q}\nonumber\\[-8pt]\\[-8pt]
&&\qquad \leq \delta_{n,s}^{(2)}
 :=
C_5 A_\cH^{4/q} B^{1/q}_\cH n^{1/2} \exp\biggl\{-
\frac{C_6}{\mathrm{f}_\infty V_{\max}^{2/s}} \biggr\}.\nonumber
\end{eqnarray}
In addition, for any $H_1\subseteq\cH$ and $H_2\subseteq\cH$
\begin{eqnarray}\label{eqmaj-1}
\bE_f
\sup_{h\in H_1}
[ \hat{r}_s(K_h)]^q &\leq&
(1+8D_s)^q
\sup_{h\in H_1} [r_s(K_h)]^q
\nonumber\\[-8pt]\\[-8pt]
&&{} +
C_7 A_\cH^2 B_\cH n^{q(s-2)/(2s)}
\exp\{-C_8 b_{n,s} \},\nonumber
\\\label{eqmaj-2}
\bE_f
\sup_{(h,\eta)\in H_1\times H_2}
[ \hat{r}_s(K_h*K_\eta)]^q
&\leq&
(1+8D_s)^q
\sup_{(h,\eta)\in H_1\times H_2} [r_s(K_h*K_\eta)]^q
\nonumber\\[-8pt]\\[-8pt]
&&{} +
C_9 A_\cH^4 B_\cH n^{q(s-2)/(2s)}
\exp\{-C_{10}b_{n,s}\},\nonumber
\end{eqnarray}
where
$b_{n,s}:=n^{4/s-1}$ if $s\in(2,4)$ and
$b_{n,s}:=[\mathrm{f}_\infty V_{\max}^{4/s}]^{-1}$ if $s\in[4,\infty)$.
The constants $C_i$, $i=2,\ldots, 10$, depend
on $L_K$, $\mathrm{k}_\infty$, $d$, $q$ and $s$ only, while $C_1$
depends also on $\mathrm{f}_\infty$.
\end{lemma}
%

\subsection{\texorpdfstring{Proofs of Theorems \protect\ref{ths<=2} and \protect\ref{ths>2}}{Proofs of Theorems 1 and 2}}
The proofs of both theorems (which we break
into several steps) follow along the same lines.

We note that in the case $s\in[2,\infty)$ the condition $f\in\bF$
implies that
$f\in\bL_s(\bR^d)$. If $s\in(1,2)$, then
by Assumptions (K1) and (K2), we have that $\bP_f\{\hat{f}_h\in
\bL_s(\bR^d)\}=1$ for any $\cX_n$-measurable vector
$h\in\cH$ and for any $n$. Hence, if $f\notin\bL_s(\bR^d)$, then
$\cR[\hat{f}_h;f]=+\infty$, $\forall h\in\cH$, and the result (i) of
Theorem~\ref{ths<=2} holds trivially. Thus, we can assume that $f \in
\bL_s(\bR^d)$ when $s\in(1,2)$.

1$^{\circ}$.
First we show that
for any $h, \eta\in\cH$
%
\begin{eqnarray}
B_{h, \eta}(f, x) &=&
B_\eta(f, x) + \int K_\eta(y-x) B_h(f,y) \,\rd y
\label{eqbias-1}
\\
&=&
B_h(f, x)+ \int K_h(y-x) B_\eta(f, y) \,\rd y.
\label{eqbias-2}
\end{eqnarray}
Indeed, by the Fubini theorem,
\begin{eqnarray*}
&& \int[K_h*K_\eta](t-x) f(t) \,\rd t\\
&&\qquad  =
\int \biggl[\int K_h(t-y) K_\eta(y-x)\,\rd y \biggr] f(t) \,\rd t
\\
&&\qquad  =
\int \biggl[\int K_h(t-y) f(t) \,\rd t - f(y) \biggr] K_\eta(y-x)
\,\rd y
+
\int K_\eta(y-x) f(y) \,\rd y
\\
&&\qquad  = \int K_\eta(y-x) f(y) \,\rd y
+ \int K_\eta(y-x) B_h(f, y) \,\rd y.
\end{eqnarray*}
Subtracting $f(x)$ from both sides of the last equality
we come to
(\ref{eqbias-1});
(\ref{eqbias-2}) follows similarly.

2$^{\circ}$.
Let $m_s(\cdot, \cdot)$ and $m_s^*(\cdot)$ be given by (\ref
{eqmajor}), and
define
%
\begin{equation}\label{eqmajorant-definition}
\delta_{n,s} := \Bigl\{\bE_f \sup_{(h, \eta)\in\cH\times\cH}
[\|\xi_{h, \eta}-\xi_{\eta}\|_s - m_{s}(h, \eta) ]_+^q \Bigr\}^{1/q}.
\end{equation}
Let $\hat{f}=\hat{f}_{\hat{h}}$ be the estimator
defined in (\ref{eqrule-0})--(\ref{eqrule}). Our first goal is to
prove that
%
\begin{equation}
\label{eqoracle-ineq-2}
\cR_s[\hat{f};f] \leq \inf_{h\in\cH}
\{ (1+3\|K\|_1 ) \cR_s[\hat{f}_{h}; f] +
3 (\bE_f [ m_s^* (h) ]^q
)^{1/q} \} + 3 \delta_{n,s}.
\end{equation}
By the triangle inequality for any $\eta\in\cH$
%
\begin{equation}\label{eqtriangle}
\|\hat{f}_{\hat{h}} - f\|_s \leq\|\hat{f}_{\hat{h}} -
\hat{f}_{\hat{h}, \eta} \|_s + \|\hat{f}_{\hat{h}, \eta} - \hat
{f}_\eta\|_s
+\|\hat{f}_\eta- f\|_s,
\end{equation}
and we are going to bound the first two terms on the right-hand side.

Define
\[
\bar{B}_{h}(f) :=
\sup_{\eta\in\cH} \biggl\| \int K_\eta(t-\cdot) B_{h}(f, t)
\,\rd t \biggr\|_s,\qquad  h\in\cH.
\]
We have for any $h\in\cH$
\begin{eqnarray*}
\hat{R}_h - m^*_s(h) &=&
\sup_{\eta\in\cH}
[ \|\hat{f}_{h, \eta}- \hat{f}_\eta\|_s - m_s(h, \eta) ]
\\
&\leq& \sup_{\eta\in\cH}
[ \|B_{h, \eta}(f, \cdot) - B_\eta(f, \cdot)\|_s +
\|\xi_{h, \eta}-\xi_\eta\|_s - m_s(h, \eta) ]
\\
&\leq& \bar{B}_h(f) +
\sup_{\eta\in\cH} [ \|\xi_{h, \eta}-\xi_\eta\|_s - m_s(h, \eta
) ]_+
=: \bar{B}_h(f) + \zeta.
\end{eqnarray*}
Here the second line is by the triangle inequality and the third line
is by
(\ref{eqbias-1}) and definition of $\bar{B}_h(f)$. Therefore for any
$h\in\cH$ one has
%
\begin{equation}\label{eqR-bound}
\hat{R}_h \leq\bar{B}_h(f) + m^*_s(h) + \zeta.
\end{equation}
By (\ref{eqbias-2}) for any $h, \eta\in\cH$
\begin{eqnarray*}
\|\hat{f}_{h, \eta}-\hat{f}_h\|_s &\leq& \|B_{h, \eta}(f, \cdot)-
B_h(f, \cdot)\|_s + \|\xi_{h, \eta}-\xi_h\|_s
\\
&\leq& \bar{B}_\eta(f) + \zeta+
\sup_{\eta\in\cH} m_s(\eta, h)
\\
&=& \bar{B}_\eta(f) + m_s^*(h) + \zeta
\leq \bar{B}_\eta(f) + \hat{R}_h + \zeta,
\end{eqnarray*}
where the last inequality is by definition of $\hat{R}_h$.
In particular, letting $h=\hat{h}$ we have
that for any $\eta\in\cH$
\begin{eqnarray}\label{eqI-1}
\|\hat{f}_{\hat{h}, \eta} - \hat{f}_{\hat{h}}\|_s
&\leq& \bar{B}_{\eta}(f) +
\hat{R}_{\hat{h}} + \zeta
\nonumber\\[-8pt]\\[-8pt]
&\leq&\bar{B}_{\eta}(f) +
\hat{R}_{\eta} + \zeta
\leq 2\bar{B}_{\eta}(f) + m^*_s(\eta) + 2\zeta,
\nonumber
\end{eqnarray}
where we have used that
$\hat{R}_{\hat{h}}\leq\hat{R}_\eta$, $\forall\eta\in\cH$
and (\ref{eqR-bound}).

Furthermore, for any $\eta\in\cH$
%
\begin{eqnarray}
\|\hat{f}_{\hat{h}, \eta} - \hat{f}_\eta\|_s &= &
\|\hat{f}_{\hat{h}, \eta} - \hat{f}_\eta\|_s - m_s(\hat{h}, \eta
) +
m_s(\hat{h}, \eta)
\nonumber
\\
&\leq& \hat{R}_{\hat{h}} + m_s^*(\eta)
\leq
\hat{R}_{\eta} + m_s^*(\eta) \leq\bar{B}_\eta(f) +
2 m_s^*(\eta) + \zeta,
\label{eqI-2}
\end{eqnarray}
where the first inequality is by definition of
$\hat{R}_h$ and $m_s^*(\cdot)$, the second inequality holds by
definition of
$\hat{h}$, and the
last inequality follows
from (\ref{eqR-bound}).

Combining (\ref{eqtriangle}), (\ref{eqI-1})
and (\ref{eqI-2}) we get for any $\eta\in\cH$ that
\begin{eqnarray*}
\|\hat{f}_{\hat{h}} - f\|_s &\leq& \|\hat{f}_{\hat{h}} - \hat
{f}_{\hat{h},
\eta}\|_s + \|\hat{f}_{\hat{h}, \eta} - \hat{f}_\eta\|_s + \|\hat
{f}_\eta-f\|_s
\\
&\leq& \|\hat{f}_\eta-f\|_s + 3\bar{B}_\eta(f) + 3m_s^*(\eta) +
3\zeta.
\end{eqnarray*}
Taking this expression to the power $q$,
computing the expectation and using the fact that $[\bE_f
|\zeta|^q]^{1/q}=\delta_{n,s}$
we obtain
%
\begin{equation}\label{eqoracle-ineq}
\mathcal{R}_s[\hat{f}; f] \leq\inf_{h\in\cH}
\{ \mathcal{R}_s[\hat{f}_{h}; f] + 3 \bar{B}_{h}(f) +
3 (\bE_f [ m_s^* (h) ]^q
)^{1/q} \} + 3 \delta_{n,s}.
\end{equation}
By the Young inequality
\[
\bar{B}_{h}(f)\leq \Bigl(\sup_{\eta\in\cH}\|K_\eta\|_1 \Bigr)
\|B_{h}(f,\cdot)\|_s=\|K\|_1 \|B_{h}(f,\cdot)\|_s.
\]
In addition [see (\ref{eq02app})--(\ref{eq03app})],
\[
\|B_{h}(f,\cdot)\|_s \leq\cR_s[\hat{f}_{h}; f]\qquad  \forall
h\in\cH.
\]
Combining this with
(\ref{eqoracle-ineq}),
we complete the proof of (\ref{eqoracle-ineq-2}).

3$^{\circ}$. Lemmas~\ref{props<=2} and~\ref{props>2} lead to an upper bound
on the
quantity $\delta_{n,s}$ given in (\ref{eqmajorant-definition}).
Indeed, by definition of $m_s(\cdot,\cdot)$ [see (\ref{eqmajor})] we have
%
\begin{eqnarray}\label{eqdelta-bound}
\delta_{n,s} & = & \Bigl\{\bE_f \sup_{(h,\eta)\in\cH\times\cH}
[\|\xi_{h,\eta}-\xi_\eta\|_s - m_s(h,\eta) ]_+^q \Bigr\}^{1/q}
\nonumber
\\
&\leq& \Bigl\{\bE_f \sup_{(h,\eta)\in\cH\times\cH}
[\|\xi_{h,\eta}\|_s - g_s(K_h*K_\eta) ]_+^q \Bigr\}^{1/q}
\\
&&{}
+
\Bigl\{\bE_f \sup_{h\in\cH}
[\|\xi_{h}\|_s - g_s(K_h) ]_+^q \Bigr\}^{1/q}
\leq \delta_{n,s}^{(1)} + \delta_{n,s}^{(2)},
\nonumber
\end{eqnarray}
where expressions for $\delta_{n,s}^{(1)}$ and $\delta_{n,s}^{(2)}$
depending on
the value of $s\in[1,\infty)$
are given in
(\ref{eqdelta-1-s<2})--(\ref{eqdelta-2-s<2}),
(\ref{eqdelta-2-1})--(\ref{eqdelta-2-2}) and
(\ref{eqlevel-1})--(\ref{eqlevel-2}).

In order to apply (\ref{eqoracle-ineq-2}) it remains to
bound
$\{\bE_f [m_s^*(h)]^q\}^{1/q}$.

4$^{\circ}$. We start with the case $s\in[1,2)$. Here, by definition,
\begin{eqnarray*}
m_s^{*}(h) &=& \sup_{\eta\in\cH} m_s(\eta, h) = g_s(K_h) +
\sup_{\eta\in\cH}
g_s(K_\eta* K_h)
\\
& =& 128 n^{1/s-1}\Bigl(\|K_h\|_s + \sup_{\eta\in\cH}\|K_h*K_\eta\|_s
\Bigr)\leq
128 [1+\|K\|_1] \mathrm{k}_\infty(nV_h)^{1/s-1}.
\end{eqnarray*}
Therefore
applying (\ref{eqoracle-ineq-2}), and taking into account
(\ref{eqdelta-bound}), (\ref{eqdelta-1-s<2}) and (\ref{eqdelta-2-s<2}),
we come to the statement (i) of Theorem~\ref{ths<=2}.

The statement (ii) of Theorem~\ref{ths<=2} dealing with the case $s=2$ follows
similarly by application of (\ref{eqoracle-ineq-2})
and (\ref{eqdelta-bound}), (\ref{eqdelta-2-1}) and (\ref{eqdelta-2-2}).
This completes the proof of Theorem~\ref{ths<=2}.

5$^{\circ}$. Now consider the case $s\in(2,\infty)$.
Because
\begin{eqnarray}\label{eqm*-BS}
m_s^*(h) &=& \sup_{\eta\in\cH} m_s(\eta, h) = g_s(K_h) +
\sup_{\eta\in\cH}
g_s(K_\eta* K_h)
\nonumber\\[-8pt]\\[-8pt]
&=& 32 \hat{r}_s(K_h) + 32 \sup_{\eta\in\cH} \hat{r}_s(K_\eta*K_h),\nonumber
\end{eqnarray}
it suffices to
bound from above $[\bE_f |\hat{r}_s(K_h)|^q]^{1/q}$ and
$[\bE_f \sup_{\eta\in\cH}|\hat{r}_s(K_h*K_\eta)|^q]^{1/q}$.
Using (\ref{eqmaj-1}) of Lemma~\ref{props>2} with
$H_1=\{h\}$ we have
\begin{eqnarray*}
[\bE_f |\hat{r}_s(K_h)|^q]^{1/q} &\leq& c_1 r_s(K_h)
+ c_2 A_\cH^{2/q}B_\cH^{1/q} n^{(s-2)/(2s)}
\exp\{-c_3 b_{n,s}\}.
\end{eqnarray*}
In addition, by the Young inequality,
\begin{eqnarray*}
\rho_s(K_h) &=& D_s n^{-1/2}\|K_h^2*f\|_{s/2}^{1/2} + n^{1/s-1}\|K_h\|_s
\\
&\leq& D_sn^{-1/2}\|K_h\|_2 \bigl\|\sqrt{f}\bigr\|_s
+ (nV_h)^{-1+1/s}\|K\|_s
\\
&\leq& D_s\mathrm{f}_\infty^{1/2}\|K\|_2 (nV_h)^{-1/2} + \|K\|_s
(nV_h)^{-1+1/s}
\leq c_4 \mathrm{f}_\infty^{1/2}(nV_h)^{-1/2};
\end{eqnarray*}
here we have used that
$
\|\sqrt{f}\|_s =(\int f^{s/2}(x) \,\rd x)^{1/s} \leq
(\mathrm{f}_\infty^{s/2-1} \int f(x)\,\rd x)^{1/s} \leq\mathrm
{f}_\infty^{1/2}.
$
Hence
\begin{eqnarray}\label{eqrss}
[\bE_f |\hat{r}_s(K_h)|^q]^{1/q} &\leq& c_5
\mathrm{f}_\infty^{1/2}(nV_h)^{-1/2}
\nonumber\\[-8pt]\\[-8pt]
&&{}
+ c_2 A_\cH^{2/q} B_\cH^{1/q} n^{(s-2)/(2s)}
\exp\{-c_3 b_{n,s}\}.
\nonumber
\end{eqnarray}

Now, applying (\ref{eqmaj-2}) with $H_1=\{h\}$ and $H_2=\cH$ we obtain
\begin{eqnarray*}
\Bigl[\bE_f \sup_{\eta\in\cH} |\hat{r}_s(K_h*K_\eta)|^q \Bigr]^{1/q} &\leq
& c_6
\sup_{\eta\in\cH}
r_s(K_h*K_\eta)
\\
&&{}
+
c_7 A_\cH^{4/q} B_\cH^{1/q} n^{(s-2)/(2s)}
\exp\{-c_{8} b_{n,s}\}.
\end{eqnarray*}
In addition, similar to the above,
\begin{eqnarray*}
\sup_{\eta\in\cH} \rho_s(K_h*K_\eta)
&\leq& \sup_{\eta\in\cH}
\bigl\{D_sn^{-1/2}\|K_h*K_\eta\|_2 \bigl\|\sqrt{f}\bigr\|_s
+ n^{-1+1/s}\|K_h*K_\eta\|_s \bigr\}
\\
& \leq&
c_8\mathrm{f}_\infty^{1/2} \sup_{\eta\in\cH} [n (V_h\vee V_\eta)]^{1/2}
\leq
c_9 \mathrm{f}_\infty^{1/2} (nV_h)^{-1/2}.
\end{eqnarray*}
Therefore the last two bounds yield
\begin{eqnarray*}
\Bigl[\bE_f \sup_{\eta\in\cH} |\hat{r}_s(K_h*K_\eta)|^q \Bigr]^{1/q} &\leq&
c_{10} \mathrm{f}_\infty^{1/2} (nV_h)^{-1/2}
\\
&&{} +
c_7 A_\cH^{4/q} B_\cH^{1/q} n^{(s-2)/(2s)}
\exp\{-c_{8} b_{n,s}\}.
\end{eqnarray*}
This along with (\ref{eqrss}) and (\ref{eqm*-BS}) results in
\begin{eqnarray*}
[\bE_f |m_s^*(K_h)|^q]^{1/q} &\leq&
c_{11} \mathrm{f}_\infty^{1/2} (nV_h)^{-1/2}
\\
&&{}
+
c_{12} A_\cH^{4/q} B_\cH^{1/q} n^{(s-2)/(2s)}
\exp\{-c_{13} b_{n,s}\}.
\end{eqnarray*}
Combining this bound with (\ref{eqlevel-1}), (\ref{eqlevel-2}) and
(\ref{eqdelta-bound}), and applying (\ref{eqoracle-ineq-2}),
we complete the proof of Theorem~\ref{ths>2}.

\subsection{\texorpdfstring{Proof of Theorem \protect\ref{thoracle}}{Proof of Theorem 3}}
Throughout the proof we denote by $c_0,c_1,\ldots,$
the positive constants depending only on the kernel $K$, the index $s$
and the
quantity $\mathrm{f}_\infty$. We divide the proof into four steps.

1$^{\circ}$. Let us prove that for any $q\geq1$ and $h\in\cH$
%
\begin{eqnarray}
\label{eq01app}
3\cR_s[\hat{f}_h; f] \geq\|B_h(f, \cdot)\|_s+\bE_f\|\xi_{h}\|_s.
\end{eqnarray}
Indeed, in view of the Jensen inequality for any $q\geq1$
%
\begin{eqnarray}
\label{eq001app}
\cR_s[\hat{f}_h; f] \geq\bE_f \|\hat{f}_h-f\|_s=\bE_f
\|B_h(f,\cdot)+\xi_{h}\|_s.
\end{eqnarray}
Denote by $\bB_p(1), 1\leq p\leq\infty,$ the unit ball in
$\bL_p (\bR^{d} )$.
By the duality argument
\begin{eqnarray*}
\bE_f \|B_h(f,\cdot)+\xi_{h}\|_s=\bE_f\sup_{\ell\in\bB
_r(1)}\int
\ell(t) [B_h(f,t)+\xi_{h}(t) ]\,\rd t,\qquad  r=\frac{s}{s-1}.
\end{eqnarray*}
Let $\ell_0\in\bB_r(1)$ be such that $\|B_h(f,\cdot)\|_s=\int
\ell_0(t)B_h(f,t)\,\rd t$;
then
%
\begin{eqnarray}
\label{eq02app}
&& \bE_f \|B_h(f,\cdot)+\xi_{h}\|_s\geq
\bE_f\int\ell_0(t) [B_h(f,t)+\xi_{h}(t) ]\,\rd t=\|B_h(f,\cdot)\|_s.
\end{eqnarray}
Here we have used that $\bE_f\xi_{h}(t)=0$, $\forall t\in\bR^d$. We
also have
by the triangle inequality
%
\begin{eqnarray}
\label{eq03app}
\bE_f \|B_h(f,\cdot)+\xi_{h}\|_s\geq\bE_f\|\xi_{h}\|_s-\|
B_h(f,\cdot)\|_s.
\end{eqnarray}
Summing up the inequalities in (\ref{eq02app}) and (\ref{eq03app}) we get
%
\begin{eqnarray}
\label{eq04app}
&& \bE_f \|B_h(f,\cdot)+\xi_{h}\|_s\geq2^{-1}\bE_f\|\xi_{h}\|_s.
\end{eqnarray}
Thus, in view of (\ref{eq02app}) and (\ref{eq04app}) for any
$\alpha\in(0,1)$
%
\begin{eqnarray}
\label{eq05app}
&& \bE_f \|B_h(f,\cdot)+\xi_{h}\|_s\geq(1-\alpha)\|B_h(f,\cdot)\|_s+
2^{-1}\alpha\bE_f\|\xi_{h}\|_s.
\end{eqnarray}
Choosing $\alpha=2/3$, we arrive to (\ref{eq01app}) in view of
(\ref{eq001app}).

In view of (\ref{eq01app}), the assertion of the theorem
will follow from the statement of Theorem~\ref{ths>2} if we show that
\[
\bE_f\|\xi_{h}\|_s\geq c_0 (nV_h)^{-1/2}.
\]

2$^{\circ}$.
Let $b>0$ be a constant to be specified,
and put $a=b^{-1}\sqrt{nV_h}.$ By duality
%
\begin{eqnarray}
\label{eq1app}
\bE_f\|\xi_{h}\|_s=\bE_f\sup_{\ell\in\bB_r(1)}\int\ell(t)\xi
_{h}(t)\,\rd t,\qquad
r=\frac{s}{s-1}.
\end{eqnarray}
Define the random event
$
\cA= \{a\xi_{h}\in\bB_2(1) \},
$
and note that if $\cA$ occurs, then by the H\"{o}lder inequality
%
\begin{eqnarray}
\label{eq2app}
a g\xi_{h}\in\bB_r(1)\qquad  \forall g\in\bB_{2r/{(2-r)}}(1).
\end{eqnarray}
Recall that $s\geq2$ implies $r\in[1,2]$, and if $r=s=2$, then we
formally put
$\frac{2r}{2-r}=\infty$.

If the event $\cA$ occurs, then $\bB_r(1) \supseteq\{ag\xi_h\dvtx g\in
\bB_{2r/{(2-r)}}(1)\}$. Therefore, by
(\ref{eq1app}) and (\ref{eq2app})
\begin{eqnarray}\label{eq3app}
\bE_f\|\xi_{h}\|_s&\geq&
a\bE_f \biggl[\bI(\cA)\sup_{g\in\bB_{2r/{(2-r)}}(1)}\int g(t)\xi
^{2}_{h}(t)\,\rd
t \biggr]
\nonumber
\\
&\geq& a\sup_{g\in\bB_{2r/{(2-r)}}(1)}\bE_f \biggl[\bI(\cA)\int
g(t)\xi^{2}_{h}(t)\,\rd t \biggr]
\nonumber\\[-8pt]\\[-8pt]
&=& a\sup_{g\in\bB_{2r/{(2-r)}}(1)}\int
g(t) [\bE_f \bI(\cA)\xi^{2}_{h}(t) ]\,\rd t =
a \|\bE_f\xi^{2}_{h}(\cdot) \bI(\cA) \|_{2s/{(s+2)}}
\nonumber\\
&\geq& a \bigl[ \|\bE_f\xi^{2}_{h}(\cdot) \|_{2s/{(s+2)}}-
\|\bE_f\xi^{2}_{h}(\cdot) \bI(\bar{\cA} ) \|_{2s/{(s+2)}}
\bigr],\nonumber
\end{eqnarray}
where $\bar{\cA}$ is the event complementary to $\cA$.

Now consider separately two cases: $s=2$ and $s>2$.

3$^{\circ}$.
If $s=2$, we get from (\ref{eq3app})
%
\begin{equation}\label{eq4app}
\bE_f\|\xi_{h}\|_2 \geq a \biggl[\int\bE_f\xi^{2}_{h}(t)\,\rd
t-\bE_f \bigl\{\|\xi_{h}\|^{2}_2 \bI\bigl(\|\xi_{h}\|_2\geq
b (nV_h)^{-1/2} \bigr) \bigr\} \biggr].
\end{equation}
Note that
%
\begin{equation}
\label{eq44app}
\bE_f\xi^{2}_{h}(t)=n^{-1}\int K^{2}_h(t-x)f(x)\,\rd x-n^{-1}
\biggl[\int
K_h(t-x)f(x)\,\rd x \biggr]^{2}
\end{equation}
and, therefore,
\[
\int
\bE_f\xi^{2}_{h}(t)\,\rd t=\frac{\|K\|^{2}_2}{nV_h}-
n^{-1}\int\biggl[\int K_h(t-x)f(x)\,\rd x \biggr]^{2}\,\rd t.
\]
The Young inequality yields
%
\begin{equation}
\label{eq444app}
\int\biggl[\int K_h(t-x)f(x)\,\rd x \biggr]^{2}\,\rd t\leq
\|K_h\|^{2}_1 \|f\|^{2}_{2}\leq\|K\|^{2}_1 \mathrm{f}_\infty.
\end{equation}
Here we have used that $f\in\bF$. Thus, in view of
$V_h\leq V_{\max}\leq1/8$ [see assumption of part (ii)
of Theorem~\ref{ths<=2}], we obtain
%
\begin{equation}
\label{eq5app}
\int\bE_f\xi^{2}_{h}(t)\,\rd
t\geq\frac{\|K\|^{2}_2}{nV_h}-\frac{\|K\|^{2}_1 \mathrm{f}_\infty
}{n}\geq c_1
(nV_h)^{-1}.
\end{equation}
It follows from Theorem 1 of GL (\citeyear{GL}) that for any $x\geq2$
%
\begin{equation}
\label{eq55app}
\bP\biggl\{\|\xi_{h}\|_2\geq\frac{x\|K\|_2}{\sqrt{nV_h}} \biggr\}\leq
e^{c_2(1-x)}
\end{equation}
and, therefore, putting $b=y\|K\|_2, y\geq2,$ we obtain
%
\begin{equation}\label{eq6app}\quad
\bE_f \biggl\{\|\xi_{h}\|^{2}_2 \bI\biggl(\|\xi_{h}\|_2\geq
\frac{y\|K\|_2}{\sqrt{nV_h}} \biggr) \biggr\}\leq
2\|K\|_2^{2}(nV_h)^{-1}\int_{y}^{\infty}xe^{c_2(1-x)}\,\rd x.
\end{equation}
Choosing $y$ sufficiently large in order to make the latter integral
less than\vspace*{-2pt}
$\frac{c_1}{4\|K\|_2^{2}}$,
we obtain from (\ref{eq4app}), (\ref{eq5app}) and (\ref{eq6app})
\[
\bE_f\|\xi_{h}\|_2 \geq c_3 (nV_h)^{-1/2}.
\]
The theorem is proved in the case $s=2$.

4$^{\circ}$.
Return now to the case $s>2.$ Note first that
\begin{eqnarray}
\label{eq7app}
\|\bE_f\xi^{2}_{h}(\cdot) \|_{{2s}/{(s+2)}}&\geq&
\biggl(\int_{B} |\bE_f\xi^{2}_{h}(t) |^{{2s}/{(s+2)}}\,\rd
t \biggr)^{{(s+2)}/{(2s)}}\nonumber\\[-8pt]\\[-8pt]
&\geq&
\nu^{{(2-s)}/{(2s)}}\int_{B}\bE_f\xi^{2}_{h}(t)\,\rd t.\nonumber
\end{eqnarray}
The last relation is obtained by the reversed H\"{o}lder inequality.
Taking into
account that $\int_{B}f(t)\,\rd t\geq\mu$, we get, using (\ref
{eq44app}) and
(\ref{eq444app}),
%
\begin{eqnarray}
\label{eq8app}
\int_{B}\bE_f\xi^{2}_{h}(t)\,\rd t\geq\frac{\mu\|K\|^{2}_2}{nV_h}-\frac{\|K\|^{2}_1 \mathrm{f}_\infty
}{n}\geq c_4\mu(nV_h)^{-1}.
\end{eqnarray}
Here we have used that $V_h\leq2^{-1}\mu\|K\|^2_2/\|K\|_1^{2}$.
On the other hand,
\[
\bE_f\xi^{2}_{h}(\cdot) \bI(\bar{\cA} )\leq
\bigl\{\bE_f [\xi_{h}(\cdot) ]^{4s/{(s+2)}} \bigr\}^{(s+2)/{(2s)}
} \{\bP(\bar{\cA}) \}^{{(s-2)}/{(2s)}}
\]
and, therefore,
\begin{eqnarray}\label{eq9app}
&&\|\bE_f\xi^{2}_{h}(\cdot) \bI(\bar{\cA} ) \|_{
{2s}/{(s+2)}}\nonumber\\[-8pt]\\[-8pt]
&&\qquad \leq
\bigl\{\bE_f \bigl( \|\xi_{h} \|_{4s/{(s+2)}} \bigr)^{4s/{(s+2)}}
\bigr\}^{(s+2)/{(2s)}} \{\bP(\bar{\cA}) \}^{(s-2)/(2s)}.\nonumber
\end{eqnarray}
We derive from Theorem 1 in GL (\citeyear{GL})
that there exists $c_5$ such that
%
\begin{equation}
\label{eq10app}
\bE_f \bigl( \|\xi_{h} \|_{4s/{(s+2)}} \bigr)^{4s/{(s+2)}}\leq c_5
(nV_h)^{-2s/{(s+2)}}.
\end{equation}
Putting $b=x\|K\|_2, x\geq2,$ we have in view of (\ref{eq55app})
\begin{eqnarray*}
\{\bP(\bar{\cA}) \}^{(s-2)/{(2s)}}\leq e^{c_2(1-x)(s-2)/{(2s)}}.
\end{eqnarray*}
It leads, together with (\ref{eq9app}) and (\ref{eq10app}), to the following
estimate:
%
\begin{eqnarray}
\label{eq11app}
\|\bE_f\xi^{2}_{h}(\cdot) \bI(\bar{\cA} ) \|_{
{2s}/{(s+2)}}\leq
c_6 (nV_h)^{-1}e^{c_2(1-x)(s-2)/{(2s)}}.
\end{eqnarray}
Finally, we obtain from (\ref{eq3app}), (\ref{eq7app}),
(\ref{eq8app}) and (\ref{eq11app})
\[
\bE_f\|\xi_{h}\|_s\geq
(x\|K\|_2)^{-1} (nV_h )^{-1/2} \bigl[c_4\mu\nu^{(2-s)/{(2s)}}-c_6
e^{c_2(1-x)(s-2)/{(2s)}} \bigr].
\]
It remains to choose $x$ sufficiently large and we come to the
assertion of the
theorem in the case $s>2$.

%
\subsection{\texorpdfstring{Proof of Theorem \protect\ref{thnikolski}}{Proof of Theorem 4}}
Let $f\in N_{s,d}(\alpha, L)$.
It is easily checked [see, e.g.,
Proposition 3
in \citeasnoun{lepski-kerk}] that bias of the estimator $\hat{f}_h$
is bounded as follows:
\[
\|B_h(f,\cdot)\|_s \leq C_1(d, s,l) L \sum_{j=1}^d h_j^{\alpha_j}.
\]
Moreover,
$ \{ \bE_f \| \xi_h\|_s^{q} \}^{1/q} \leq
C_2(nV_h)^{-\gamma_s}$.
If we
set the ``oracle bandwidth'' $h^*:=(h_1^*,\ldots, h_d^*)$ so that
\[
[ h_j^*]^{\alpha_j}:=
\biggl[\frac{C_2}{C_1} \biggr]^{\bar{\alpha}/(\gamma_s+\bar{\alpha})}
L^{-\bar{\alpha}/(\gamma_s+\bar{\alpha})}
n^{-\gamma_s\bar{\alpha}/(\gamma_s+\bar{\alpha})},\qquad  j=1,\ldots, d,
\]
then $h^*\in\cH$ and $\hat{f}_{h^*}\in\cF(\cH)$ for large enough $n$.
Hence, for any $f\in N_{s,d}(\alpha, L)$ we have that
$\cR_s[\hat{f}_{h^*}; f] \leq C_3 \varphi_{n,s}(\bar{\alpha})$.
Then we apply
oracle
inequalities of Theorems~\ref{ths<=2} and~\ref{ths>2}.
Observe that by choice of constant $\kappa_2$ in definition of
$h^{\max}$ we guarantee that the
remainder terms
are negligibly small as $n\to\infty$
in comparison with the first
terms in
(\ref{eqor-2}) and (\ref{eqor-3}). This fact leads to the statement
of the
theorem.

\subsection{\texorpdfstring{Proof of Theorem \protect\ref{thcompact}}{Proof of Theorem 5}}
First we note that it suffices to prove the theorem only for $s\geq2$.
Indeed, since $\operatorname{supp}(f)\subseteq Q$, one has
$\operatorname{supp} (\hat{f}_h )\subseteq Q^{\prime}$
for any $\cX_n$-measurable random vector $h\in\cH$, where, in view of
the
assumptions imposed on the kernel $K$,
\[
Q^{\prime}= \{y\in\bR^d: |y_i-x_i|\leq1/2, i=1,\ldots, d, x\in Q \}.
\]
Here we have also used that $h^{\max}\in(0,1]^d$.
Thus, for any density $f$ and any $\cX_n$-measurable random vector
$h\in\cH$
\[
\operatorname{supp} (\hat{f}_h-f )\subseteq Q^{\prime}
\]
and, therefore, in view of
H\"{o}lder inequality for any $s\in[1,2)$
\[
\|\hat{f}_h-f \|_s\leq[\operatorname{mes} \{Q^{\prime} \}]^{
(2-s)/{(2s)}} \|\hat{f}_h-f \|_2.
\]
We conclude that for any $s\in[1,2)$ the estimation problem in the
$\bL_s$-norm
can be reduced to the estimation problem in the $\bL_2$-norm.

Let $f\in W_{p,d}(\alpha,L,Q)$ and $s\geq2$.
The standard computation
(by the generalized Minkowski inequality and by the H\"{o}lder
inequality along
with the fact that $f$ is compactly supported) yields the following
bound on the
$\bL_s$-norm of the bias of $\hat{f}_h$:
\begin{eqnarray*}
\|B_h(f,\cdot)\|_s &\leq&
C_1(d,s,l) L [\operatorname{mes}\{Q\}]^{(p-s)/{(sp)}} \sum_{j=1}^d
h_j^{\alpha_j}.
\end{eqnarray*}
Moreover,
$ \{ \bE_f \| \xi_h\|_s^{q} \}^{1/q} \leq
C_2(nV_h)^{-1/2}$.
If we
set the ``oracle bandwidth'' $h^*:=(h_1^*,\ldots, h_d^*)$ so that
%
\begin{eqnarray} 
[ h_j^*]^{\alpha_j}:=
\biggl[\frac{C_2}{C_1} \biggr]^{2\bar{\alpha}/(2\bar{\alpha}+1)}
\bigl(L [\operatorname{mes}\{Q\} ]^{(p-s)/{(sp)}} \bigr)^{-2\bar{\alpha
}/(1+2\bar{\alpha})}
n^{-\bar{\alpha}/(2\bar{\alpha}+1)},\nonumber \\
\eqntext{j=1,\ldots, d,}
\end{eqnarray}
then $h^*\in\cH$ and $\hat{f}_{h^*}\in\cF(\cH)$ for large enough $n$.
Then the result follows by application of Theorems~\ref{ths<=2}(ii)
and~\ref{ths>2}.

\begin{appendix}
\section*{Appendix}\label{app}
Proofs of Lemmas~\ref{props<=2} and~\ref{props>2}
follow directly from general uniform bounds on norms of empirical
processes established in GL (\citeyear{GL}).
In our proofs below we use notation and terminology of this
paper.

\begin{pf*}{Proof of Lemma~\ref{props<=2}}
The statement is a direct consequence of Theorem~4 of Section 3.3 in
GL (\citeyear{GL}).

To apply this theorem
one should verify Assumptions (W1), (W4) and (L)
for the following classes of weights
$\cW^{(1)}=\{w=n^{-1}K_h\dvtx h\in\cH\}$ and $\cW^{(2)}=\{
w=n^{-1}(K_h*K_\eta)\dvtx
(h,\eta)\in\cH\times\cH\}$.
The sets $\cW^{(1)}$ and $\cW^{(2)}$ are considered as images of $\cH
$ and
$\cH\times\cH$
under transformations $h\mapsto n^{-1}K_h$ and $(h, \eta)\mapsto
n^{-1}(K_h*K_\eta)$, respectively.
The sets $\cH$ and $\cH\times\cH$ are equipped with the distances
\begin{eqnarray*}
\mathrm{d}_1(h, h^\prime)&\hspace*{3pt}=&c_1 \max_{i=1,\ldots,d} \ln\biggl(\frac
{h_i\vee
h_i^\prime}{h_i\wedge h_i^\prime}
\biggr),\\
\mathrm{d}_2[(h,h^\prime), (\eta,\eta^\prime)]&:=& c_2 \{
\mathrm{d}_1(h, h^\prime)\vee\mathrm{d}_1(\eta,\eta^\prime)\},
\end{eqnarray*}
where $c_1$ and $c_2$ are appropriate constants depending on
$\mathrm{k}_\infty$,
$L_K$ and $d$ only [see formulas (9.1) and (9.2) in GL (\citeyear{GL})].
With this notation Lemma 9 of GL (\citeyear{GL}) shows that Assumption
(L) holds
for both $\cW^{(1)}$ and
$\cW^{(2)}$. Moreover,
Assumption~(W1) holds trivially for both $\cW^{(1)}$ and $\cW^{(2)}$ with
$\mu_*=V_{\max}$ and $\mu_*=2^d V_{\max}$, respectively.
Moreover, Assumption (W4) for both $\cW^{(1)}$ and $\cW^{(2)}$
follows from
formula (9.8) in GL (\citeyear{GL}).
Thus all conditions of Theorem 4 are fulfilled.

(i)
We apply this theorem with $z=1$ and $\epsilon=1$.
We need to evaluate the constant $T_{3,\epsilon}$
for $\cW^{(1)}$ and $\cW^{(2)}$. If $N_{\cH, \mathrm{d}_1}(\epsilon
)$ denotes
the minimal number of balls in the metric $\mathrm{d}_1$ needed to
cover $\cH$,
then formula (9.8)
from GL (\citeyear{GL})
shows that
$N_{\cH, \mathrm{d_1}}(1/8) \leq c_3 A_{\cH}$,
where $c_3$ depends on $d$ only. Similarly, $N_{\cH\times\cH,
\mathrm{d}_2}(1/8)\leq c_4 A_{\cH}^2$.
In addition, for
\[
L_{\cH, \mathrm{d}_1}(\epsilon) :=
\sum_{k=1}^\infty\exp\{2\ln N_{\cH,\mathrm{d}_1}(\epsilon2^{-k})
-(9/16)2^k k^{-2} \}
\]
we have $L_{\cH, \mathrm{d}_1}(1) \leq c_5 A_\cH$.
Similarly, $L_{\cH\times\cH, \mathrm{d}_2}(1) \leq c_6 A^2_\cH$.
Combining these bounds we come to the statement (i).

(ii) The second statement follows exactly in the same way from the above
considerations.
Theorem 4 of GL (\citeyear{GL}) is again applied with $z=1$ and $\epsilon=1$.
\end{pf*}

\begin{pf*}{Proof of Lemma~\ref{props>2}}
The proof is by
application of Theorem 7 from GL (\citeyear{GL}).
We need to calculate several
quantities.

We start with the class $\cW^{(1)}$.
Here for $\vartheta_0^{(1)}=10D_s\mathrm{f}_\infty(L_K \sqrt{d})^{d/2}$
we have
\begin{eqnarray*}
C^*_{\xi,1}(y)&=&1+2\vartheta_0^{(1)} \bigl\{\sqrt{y} \bigl(V_{\max}^{1/s}+n^{
-1/(2s)} \bigr)+
yn^{-1/s} \bigr\}
\\
&\leq& 1+2\vartheta_0^{(1)} \bigl\{2\sqrt{y}V_{\max}^{1/s}+yV_{\max
}^{2/s} \bigr\},
\end{eqnarray*}
where we have used that $V_{\max}\geq1/\sqrt{n}$.
If we set $y=\bar{y}:=[4V_{\max}^{2/s} (\vartheta_0^{(1)}\vee
1)]^{-1}$, then
$C^*_{\xi,1}(\bar{y})\leq4$.
We apply Theorem 7 with $\epsilon=1$ and $y=\bar{y}$. Condition
$nV_{\min}> C_1=[256 D_s^2]^{(s\wedge4)/(s\wedge4 -2)}$ implies that
\[
\bar{u}_1(\gamma)=4 \bigl[1- 8 D_s (nV_{\min})^{1/(s\wedge
4)-1/2} \bigr]^{-1}\leq8.
\]
Moreover, we note that condition $\bar{y}\leq y_*^{(1)}$
follows from definition of $\bar{y}$ and $n\geq C_2$.
In addition,
$\tilde{T}_{1,\epsilon}^{(1)}\leq cA_\cH^2 B_\cH$.
These facts imply (\ref{eqlevel-1}) and (\ref{eqmaj-1}).

The bounds (\ref{eqlevel-2}) and (\ref{eqmaj-2}) for $\cW^{(2)}$
follow from
similar computations.
\end{pf*}
\end{appendix}

\section*{Acknowledgments} The authors thank two anonymous referees for
useful comments and suggestions.


%

\printaddresses


\begin{thebibliography}{27}

\bibitem[\protect\citeauthoryear{Birg{\'e}}{2008}]{birge}
%
\begin{bmisc}[auto:STB|2011-03-03|12:04:44]
\bauthor{\bsnm{Birg{\'e}},~\bfnm{L.}\binits{L.}}
(\byear{2008}).
\bhowpublished{Model selection for density estimation with $\mathbb
{L}_2$-loss.
Available at
\href{http://arxiv.org/abs/arXiv:0808.1416v2}{arXiv:0808.1416v2}}.
\end{bmisc}
%
\endbibitem

\bibitem[\protect\citeauthoryear{Bretagnolle and Huber}{1979}]{bretagnolle}
%
\begin{barticle}[mr]
\bauthor{\bsnm{Bretagnolle},~\bfnm{J.}\binits{J.}} \AND
\bauthor{\bsnm{Huber},~\bfnm{C.}\binits{C.}}
(\byear{1979}).
\btitle{Estimation des densit\'es: Risque minimax}.
\bjournal{Z.~Wahrsch. Verw. Gebiete}
\bvolume{47}
\bpages{119--137}.
\bid{doi={10.1007/BF00535278}, issn={0044-3719}, mr={0523165}}
\end{barticle}
%
\endbibitem

\bibitem[\protect\citeauthoryear{Devroye and Gy{\"
o}rfi}{1985}]{devroye-gyorfi}
%
\begin{bbook}[mr]
\bauthor{\bsnm{Devroye},~\bfnm{Luc}\binits{L.}} \AND
\bauthor{\bsnm{Gy{\"o}rfi},~\bfnm{L{\'a}szl{\'o}}\binits{L.}}
(\byear{1985}).
\btitle{Nonparametric Density Estimation: The $L{\sb{1}}$ View}.
\bpublisher{Wiley}, \baddress{New York}.
\bid{mr={0780746}}
\end{bbook}
%
\endbibitem

\bibitem[\protect\citeauthoryear{Devroye and Lugosi}{1996}]{dev-lug96}
%
\begin{barticle}[mr]
\bauthor{\bsnm{Devroye},~\bfnm{Luc}\binits{L.}} \AND
\bauthor{\bsnm{Lugosi},~\bfnm{G{\'a}bor}\binits{G.}}
(\byear{1996}).
\btitle{A universally acceptable smoothing factor for kernel density
estimates}.
\bjournal{Ann. Statist.}
\bvolume{24}
\bpages{2499--2512}.
\bid{doi={10.1214/aos/1032181164}, issn={0090-5364}, mr={1425963}}
\end{barticle}
%
\endbibitem

\bibitem[\protect\citeauthoryear{Devroye and Lugosi}{1997}]{dev-lug97}
%
\begin{barticle}[mr]
\bauthor{\bsnm{Devroye},~\bfnm{Luc}\binits{L.}} \AND
\bauthor{\bsnm{Lugosi},~\bfnm{G{\'a}bor}\binits{G.}}
(\byear{1997}).
\btitle{Nonasymptotic universal smoothing factors, kernel complexity and
{Y}atracos classes}.
\bjournal{Ann. Statist.}
\bvolume{25}
\bpages{2626--2637}.
\bid{doi={10.1214/aos/1030741088}, issn={0090-5364}, mr={1604428}}
\end{barticle}
%
\endbibitem

\bibitem[\protect\citeauthoryear{Devroye and Lugosi}{2001}]{devroye-lugosi}
%
\begin{bbook}[mr]
\bauthor{\bsnm{Devroye},~\bfnm{Luc}\binits{L.}} \AND
\bauthor{\bsnm{Lugosi},~\bfnm{G{\'a}bor}\binits{G.}}
(\byear{2001}).
\btitle{Combinatorial Methods in Density Estimation}.
\bpublisher{Springer}, \baddress{New York}.
\bid{mr={1843146}}
\end{bbook}
%
\endbibitem

\bibitem[\protect\citeauthoryear{Donoho et~al.}{1996}]{Donoho}
%
\begin{barticle}[mr]
\bauthor{\bsnm{Donoho},~\bfnm{David~L.}\binits{D.~L.}},
\bauthor{\bsnm{Johnstone},~\bfnm{Iain~M.}\binits{I.~M.}},
\bauthor{\bsnm{Kerkyacharian},~\bfnm{G{\'e}rard}\binits{G.}} \AND
\bauthor{\bsnm{Picard},~\bfnm{Dominique}\binits{D.}}
(\byear{1996}).
\btitle{Density estimation by wavelet thresholding}.
\bjournal{Ann. Statist.}
\bvolume{24}
\bpages{508--539}.
\bid{doi={10.1214/aos/1032894451}, issn={0090-5364}, mr={1394974}}
\end{barticle}
%
\endbibitem

\bibitem[\protect\citeauthoryear{Goldenshluger and Lepski}{2008}]{GL-1}
%
\begin{barticle}[mr]
\bauthor{\bsnm{Goldenshluger},~\bfnm{Alexander}\binits{A.}} \AND
\bauthor{\bsnm{Lepski},~\bfnm{Oleg}\binits{O.}}
(\byear{2008}).
\btitle{Universal pointwise selection rule in multivariate function
estimation}.
\bjournal{Bernoulli}
\bvolume{14}
\bpages{1150--1190}.
\bid{doi={10.3150/08-BEJ144}, issn={1350-7265}, mr={2543590}}
\end{barticle}
%
\endbibitem

\bibitem[\protect\citeauthoryear{Goldenshluger and Lepski}{2009}]{GL-2}
%
\begin{barticle}[mr]
\bauthor{\bsnm{Goldenshluger},~\bfnm{Alexander}\binits{A.}} \AND
\bauthor{\bsnm{Lepski},~\bfnm{Oleg}\binits{O.}}
(\byear{2009}).
\btitle{Structural adaptation via {$\Bbb L\sb p$}-norm oracle inequalities}.
\bjournal{Probab. Theory Related Fields}
\bvolume{143}
\bpages{41--71}.
\bid{doi={10.1007/s00440-007-0119-5}, issn={0178-8051}, mr={2449122}}
\end{barticle}
%
\endbibitem

\bibitem[\protect\citeauthoryear{Goldenshluger and Lepski}{2011}]{GL}
%
\begin{bmisc}[auto:STB|2011-03-03|12:04:44]
\bauthor{\bsnm{Goldenshluger},~\bfnm{A.}\binits{A.}} \AND
\bauthor{\bsnm{Lepski},~\bfnm{O.}\binits{O.}}
(\byear{2011}).
\bhowpublished{Uniform bounds for norms of sums of independent random
functions. \textit{Ann. Probab.} To appear. Available at
\href{http://arxiv.org/abs/arXiv:0904.1950v2}{arXiv:0904.1950v2}.}
\end{bmisc}
%
\endbibitem

\bibitem[\protect\citeauthoryear{Hasminskii and Ibragimov}{1990}]{Has-Ibr}
%
\begin{barticle}[mr]
\bauthor{\bsnm{Hasminskii},~\bfnm{Rafael}\binits{R.}} \AND
\bauthor{\bsnm{Ibragimov},~\bfnm{Ildar}\binits{I.}}
(\byear{1990}).
\btitle{On density estimation in the view of {K}olmogorov's ideas in
approximation theory}.
\bjournal{Ann. Statist.}
\bvolume{18}
\bpages{999--1010}.
\bid{doi={10.1214/aos/1176347736}, issn={0090-5364}, mr={1062695}}
\end{barticle}
%
\endbibitem

\bibitem[\protect\citeauthoryear{Ibragimov and
Has'minski{\u\i}}{1980}]{Ibr-Has1}
%
\begin{barticle}[mr]
\bauthor{\bsnm{Ibragimov},~\bfnm{I.~A.}\binits{I.~A.}} \AND
\bauthor{\bsnm{Has'minski{\u\i}},~\bfnm{R.~Z.}\binits{R.~Z.}}
(\byear{1980}).
\btitle{An estimate of the density of a distribution}.
\bjournal{Zap. Nauchn. Sem. Leningrad. Otdel. Mat. Inst. Steklov. (LOMI)}
\bvolume{98}
\bpages{61--85}.
\bid{issn={0207-6772}, mr={0591862}}
\end{barticle}
%
\endbibitem

\bibitem[\protect\citeauthoryear{Ibragimov and
Khas'minski{\u\i}}{1981}]{Ibr-Has2}
%
\begin{barticle}[mr]
\bauthor{\bsnm{Ibragimov},~\bfnm{I.~A.}\binits{I.~A.}} \AND
\bauthor{\bsnm{Khas'minski{\u\i}},~\bfnm{R.~Z.}\binits{R.~Z.}}
(\byear{1981}).
\btitle{More on estimation of the density of a distribution}.
\bjournal{Zap. Nauchn. Sem. Leningrad. Otdel. Mat. Inst. Steklov. (LOMI)}
\bvolume{108}
\bpages{72--88}.
\bid{issn={0206-8540}, mr={0629401}}
\end{barticle}
%
\endbibitem

\bibitem[\protect\citeauthoryear{Jennrich}{1969}]{jennrich}
%
\begin{barticle}[mr]
\bauthor{\bsnm{Jennrich},~\bfnm{Robert~I.}\binits{R.~I.}}
(\byear{1969}).
\btitle{Asymptotic properties of non-linear least squares estimators}.
\bjournal{Ann. Math. Statist.}
\bvolume{40}
\bpages{633--643}.
\bid{issn={0003-4851}, mr={0238419}}
\end{barticle}
%
\endbibitem

\bibitem[\protect\citeauthoryear{Johnson, Schechtman and
Zinn}{1985}]{Johnson}
%
\begin{barticle}[mr]
\bauthor{\bsnm{Johnson},~\bfnm{W.~B.}\binits{W.~B.}},
\bauthor{\bsnm{Schechtman},~\bfnm{G.}\binits{G.}} \AND
\bauthor{\bsnm{Zinn},~\bfnm{J.}\binits{J.}}
(\byear{1985}).
\btitle{Best constants in moment inequalities for linear combinations of
independent and exchangeable random variables}.
\bjournal{Ann. Probab.}
\bvolume{13}
\bpages{234--253}.
\bid{issn={0091-1798}, mr={0770640}}
\end{barticle}
%
\endbibitem

\bibitem[\protect\citeauthoryear{Juditsky and
Lambert-Lacroix}{2004}]{Juditsky}
%
\begin{barticle}[mr]
\bauthor{\bsnm{Juditsky},~\bfnm{Anatoli}\binits{A.}} \AND
\bauthor{\bsnm{Lambert-Lacroix},~\bfnm{Sophie}\binits{S.}}
(\byear{2004}).
\btitle{On minimax density estimation on {$\Bbb R$}}.
\bjournal{Bernoulli}
\bvolume{10}
\bpages{187--220}.
\bid{doi={10.3150/bj/1082380217}, issn={1350-7265}, mr={2046772}}
\end{barticle}
%
\endbibitem

\bibitem[\protect\citeauthoryear{Kerkyacharian, Lepski and
Picard}{2001}]{lepski-kerk}
%
\begin{barticle}[mr]
\bauthor{\bsnm{Kerkyacharian},~\bfnm{G{\'e}rard}\binits{G.}},
\bauthor{\bsnm{Lepski},~\bfnm{Oleg}\binits{O.}} \AND
\bauthor{\bsnm{Picard},~\bfnm{Dominique}\binits{D.}}
(\byear{2001}).
\btitle{Nonlinear estimation in anisotropic multi-index denoising}.
\bjournal{Probab. Theory Related Fields}
\bvolume{121}
\bpages{137--170}.
\bid{doi={10.1007/PL00008800}, issn={0178-8051}, mr={1863916}}
\end{barticle}
%
\endbibitem

\bibitem[\protect\citeauthoryear{Kerkyacharian, Picard and
Tribouley}{1996}]{kerk}
%
\begin{barticle}[mr]
\bauthor{\bsnm{Kerkyacharian},~\bfnm{G{\'e}rard}\binits{G.}},
\bauthor{\bsnm{Picard},~\bfnm{Dominique}\binits{D.}} \AND
\bauthor{\bsnm{Tribouley},~\bfnm{Karine}\binits{K.}}
(\byear{1996}).
\btitle{{$L\sp p$} adaptive density estimation}.
\bjournal{Bernoulli}
\bvolume{2}
\bpages{229--247}.
\bid{doi={10.2307/3318521}, issn={1350-7265}, mr={1416864}}
\end{barticle}
%
\endbibitem

\bibitem[\protect\citeauthoryear{Mason}{2009}]{mason}
%
\begin{barticle}[auto:STB|2011-03-03|12:04:44]
\bauthor{\bsnm{Mason},~\bfnm{D.~M.}\binits{D.~M.}}
(\byear{2009}).
\btitle{Risk bounds for kernel density estimators}.
\bjournal{Zap. Nauchn. Sem. Leningrad. Otdel. Mat. Inst. Steklov. (LOMI)}
\bvolume{363}
\bpages{66--104}.
\bnote{Available at
\texttt{\href{http://www.pdmi.ras.ru/znsl/}{http://}
\href{http://www.pdmi.ras.ru/znsl/}{www.pdmi.ras.ru/znsl/}}}.
\end{barticle}
%
\endbibitem

\bibitem[\protect\citeauthoryear{Massart}{2007}]{massart}
%
\begin{bbook}[mr]
\bauthor{\bsnm{Massart},~\bfnm{Pascal}\binits{P.}}
(\byear{2007}).
\btitle{Concentration Inequalities and Model Selection}.
\bseries{Lecture Notes in Math.}
\bvolume{1896}.
\bpublisher{Springer}, \baddress{Berlin}.
\bid{mr={2319879}}
\end{bbook}
%
\endbibitem

\bibitem[\protect\citeauthoryear{Nikol'ski{\u\i}}{1969}]{nikol}
%
\begin{bbook}[mr]
\bauthor{\bsnm{Nikol'ski{\u\i}},~\bfnm{S.~M.}\binits{S.~M.}}
(\byear{1969}).
\btitle{Priblizhenie Funktsii Mnogikh Peremennykh i Teoremy Vlozheniya}.
\bpublisher{Nauka}, \baddress{Moscow}.
\bid{mr={0310616}}
\bptnote{check related}%
\end{bbook}
%
\endbibitem

\bibitem[\protect\citeauthoryear{Parzen}{1962}]{parzen}
%
\begin{barticle}[mr]
\bauthor{\bsnm{Parzen},~\bfnm{Emanuel}\binits{E.}}
(\byear{1962}).
\btitle{On estimation of a probability density function and mode}.
\bjournal{Ann. Math. Statist.}
\bvolume{33}
\bpages{1065--1076}.
\bid{issn={0003-4851}, mr={0143282}}
\end{barticle}
%
\endbibitem

\bibitem[\protect\citeauthoryear{Rigollet and
Tsybakov}{2007}]{rigollet-tsybakov}
%
\begin{barticle}[mr]
\bauthor{\bsnm{Rigollet},~\bfnm{Ph.}\binits{P.}} \AND
\bauthor{\bsnm{Tsybakov},~\bfnm{A.~B.}\binits{A.~B.}}
(\byear{2007}).
\btitle{Linear and convex aggregation of density estimators}.
\bjournal{Math. Methods Statist.}
\bvolume{16}
\bpages{260--280}.
\bid{doi={10.3103/S1066530707030052}, issn={1066-5307}, mr={2356821}}
\end{barticle}
%
\endbibitem

\bibitem[\protect\citeauthoryear{Rosenblatt}{1956}]{rosenblatt}
%
\begin{barticle}[mr]
\bauthor{\bsnm{Rosenblatt},~\bfnm{Murray}\binits{M.}}
(\byear{1956}).
\btitle{Remarks on some nonparametric estimates of a density function}.
\bjournal{Ann. Math. Statist.}
\bvolume{27}
\bpages{832--837}.
\bid{issn={0003-4851}, mr={0079873}}
\end{barticle}
%
\endbibitem

\bibitem[\protect\citeauthoryear{Samarov and Tsybakov}{2007}]{samarov}
%
\begin{bincollection}[mr]
\bauthor{\bsnm{Samarov},~\bfnm{Alexander}\binits{A.}} \AND
\bauthor{\bsnm{Tsybakov},~\bfnm{Alexandre}\binits{A.}}
(\byear{2007}).
\btitle{Aggregation of density estimators and dimension reduction}.
In \bbooktitle{Advances in Statistical Modeling and Inference}
(\beditor{V. Nair}, ed.).
\bseries{Ser. Biostat.}
\bvolume{3}
\bpages{233--251}.
\bpublisher{World Scientific}, \baddress{Hackensack, NJ}.
\bid{mr={2416118}}\
\end{bincollection}
%
\endbibitem

\bibitem[\protect\citeauthoryear{Silverman}{1986}]{silverman}
%
\begin{bbook}[mr]
\bauthor{\bsnm{Silverman},~\bfnm{B.~W.}\binits{B.~W.}}
(\byear{1986}).
\btitle{Density Estimation for Statistics and Data Analysis}.
\bpublisher{Chapman \& Hall}, \baddress{London}.
\bid{mr={0848134}}
\end{bbook}
%
\endbibitem

\end{thebibliography}
\end{document}